\documentclass[11pt]{amsart}

\usepackage[square,compress,comma, numbers,sort]{natbib}
\usepackage[colorlinks=true, citecolor=red, linkcolor=blue]{hyperref}
\usepackage{amsfonts,mathtools}

\allowdisplaybreaks[4]

\usepackage{amsmath}
\usepackage{bbm}
\usepackage{amssymb}
\usepackage{mathtools}

\usepackage{color}

\definecolor{c20}{rgb}{0.,0.7,0.}
\definecolor{c30}{rgb}{0.,0.,1.}
\definecolor{c40}{rgb}{1,0.1,0.7}
\definecolor{c50}{rgb}{1,0,0}
\definecolor{c60}{rgb}{1,0.9,0.1}

\def\x{\vk{x}}

\newcommand{\ve}{\varepsilon}

\newcommand{\E}[1]{\mathbb{E}\left\{ #1\right\}}

\newcommand{\pk}[1]{\mathbb{P} \left\{ #1 \right \} }

\newcommand{\R}{\mathbb{R}}

\newcommand{\BQN}{\begin{eqnarray}}
	\newcommand{\EQN}{\end{eqnarray}}
\newcommand{\BQNY}{\begin{eqnarray*}}
	\newcommand{\EQNY}{\end{eqnarray*}}

\newcommand{\BS}{\begin{sat}}
	\newcommand{\ES}{\end{sat}}
\newcommand{\BT}{\begin{theo}}
	\newcommand{\ET}{\end{theo}}
\newcommand{\BK}{\begin{korr}}
	\newcommand{\EK}{\end{korr}}

\newcommand{\BD}{\begin{de}}
	\newcommand{\ED}{\end{de}}

\newcommand{\BIT}{\begin{itemize}}
	\newcommand{\EIT}{\end{itemize}}
\newcommand{\BDI}{\begin{description}}
	\newcommand{\EDI}{\end{description}}

\newcommand{\BRM}{\begin{remarks}}
	\newcommand{\ERM}{\end{remarks}}

\newcommand{\BEL}{\begin{lem}}
	\newcommand{\EEL}{\end{lem}}

\newtheorem{theo}{Theorem}[section]
\newtheorem{sat}[theo]{Proposition}
\newtheorem{de}[theo]{Definition}
\newtheorem{lem}[theo]{Lemma}

\newtheorem{example}[theo]{Example}
\newtheorem{korr}[theo]{Corollary}
\newtheorem{remark}[theo]{Remark}
\newtheorem{remarks}[theo]{Remarks}

\newcommand{\COM}[1]{}

\def\ve{\varepsilon}

\def\IF{\infty}

\newcommand{\kb}[1]{\boldsymbol{#1}}
\newcommand{\vk}[1]{\kb{#1}}

\def\ve{\varepsilon}

\def\IF{\infty}

\begin{document}
	
	\title{Ruin Probability Approximation for Bidimensional Brownian Risk Model with Tax}
	
	\author{Timofei Shashkov}
	\address{Timofei, Department of Actuarial Science, 
		University of Lausanne,\\
		UNIL-Dorigny, 1015 Lausanne, Switzerland
	}
	\email{timofei98shashkov@gmail.com}
	
	\bigskip
	
	\date{\today}
	\maketitle

	{\bf Abstract:} 
		Let $\vk B(t)=(B_1(t), B_2(t))$, $t\geq 0$ be a two-dimensional Brownian motion with independent components and define the $\vk\gamma$-reflected process 
		$$\vk X(t)=(X_1(t),X_2(t))=\left(B_1(t)-c_1t-\gamma_1\inf_{s_1\in[0,t]}(B_1(s_1)-c_1s_1),B_2(t)-c_2t-\gamma_2\inf_{s_2\in[0,t]}(B_2(s_2)-c_2s_2)\right),$$ with given finite constants $c_1,c_2$ and $\gamma_1,\gamma_2\in[0,2)$. The goal of this paper is to derive the asymptotics of the ruin probability
		$$\pk{\exists_{t\in[0,T]}: X_1(t)>u,X_2(t)>au}$$
		as $u\to\infty$ for $T>0$.

	

	\section{Introduction and main result}
	The classical Brownian risk model (BRM) of a given insurance portfolio plays a significant role in risk theory (see \cite{diff}). The corresponding risk process $R$ is defined by
	$$R(t)= u+ ct - B(t), \ \ \ \ \ t\geq 0,$$
	with \(B\) being a Brownian motion, $u>0$ the initial capital and $c>0$ the premium rate. \\
	For a finite time horizon $[0,T],T>0$ the classical ruin probability is defined by (see e.g., \cite{DeM15})
	\begin{align*}
		\widetilde{ \psi}_T(u)\coloneqq\pk{ \inf_{t\in [0,T]} R(t) < 0}= \pk{\sup_{t\in [0,T]} (B(t)- ct)> u},\ u>0.
	\end{align*}
	According to \cite{Krzyspp} for $T<\IF$
	\begin{align*} 
		\widetilde{ \psi}_T(u) =\Phi\left(-\frac u{\sqrt{T}} -c\sqrt{T}\right)+
		e^{-2cu}\Phi\left(- \frac u{\sqrt{T}} +c\sqrt{T}\right),
	\end{align*}
	where $\Phi$ is the distribution function of a standard normal random variable. In addition, if $T=\IF$, then
	$$\widetilde{ \psi}_\IF(u)\coloneqq\pk{ \inf_{t\geq 0} R(t) < 0}=e^{-2cu}.$$
	An extension of the Brownian risk model, which was motivated in \cite{Alb}, accounts for the non-zero tax rate \(\gamma\in[0,1]\). The corresponding risk process \(X\) is given by
	\[
	X(t)=B(t)-ct-\gamma\inf_{s\in[0,t]}(B(s)-cs), \ \ \ \ \ t\geq 0.
	\]
	This process is also known in the literature under the name of \(\gamma\)-reflected Brownian motion, see \cite{Gamma1, Gamma_pass, Ontheref, Extref, Grisha}.
	
	If $\gamma=1$, then the process $X$ is also known in queuing theory as a workload process or a queue length process and it has been studied in several contributions (see, e.g. \cite{G1, Whitt, AwGly}).
	
	
	
	Unlike the classical Brownian model for this type of process it is not possible to calculate this ruin probability explicitly. Therefore it can either be simulated or approximated for large values of the initial capital $u$. According to \cite{Gamma1}, if \(\psi_{\gamma,T}(u)=\pk{\exists_{t\in[0,T]}: X(t)>u}
	\) for \(0<T<\IF\), then
	$$\psi_{\gamma,T}(u)\sim \frac{4}{2-\gamma}\psi\left(\frac{u+cT}{\sqrt{T}}\right)\sim \frac{4}{2-\gamma}\pk{B(T)>u+cT},\ u\to\IF,$$
	where $\psi(x)=1/(\sqrt{2\pi}x)\exp(-x^2/2)$. In this paper $\sim$ means asymptotic equivalence as \(u\to\IF\).
	
	
	In many applications, vector-valued processes are frequently used, making their analysis particularly important. 
	For the previously discussed \(\gamma\)-reflected risk model, which is a benchmark model in risk theory, it is hence natural to consider its vector-valued generalization. In this manuscript under certain assumptions we define and study the bivariate Brownian risk model with tax. To this end, consider two independent one-dimensional Brownian motions, $B_1$ and $B_2$. The model is then defined by  \(\vk X(t)=(X_1(t),X_2(t))\) with
	\[
	X_i(t)=B_i(t)-c_it-\gamma_i\inf_{s_i\in [0,t]}(B_i(s_i)-c_is_i), \ \ \ \ t\geq 0
	\]
	for \(i=1,2\), where \(\gamma_1,\gamma_2\) are tax rates, \(c_1,c_2\) are premium rates. Even though in practice it is natural to assume \(c_1,c_2>0\) and \(\gamma_1,\gamma_2\in[0,1]\), methods which are used in the manuscript to obtain main results are applicable for all \(c_1,c_2\in\R\) and \(\gamma_1,\gamma_2\in[0,2)\). 
	
	The main quantities of interest related to the vector-valued process \(\vk X\) in the context of risk theory are the finite and infinite horizon ruin events. In this paper we study the simultaneous ruin probability in the case of a finite time frame, defined for \(T>0\) and \(\vk u=(u_1,u_2)\) by the formula
	\[
	\psi_{\gamma, T}(\vk u)=\pk{\exists_{t\in [0,T]}: X_1(t)>u_1, X_2(t)>u_2}.
	\]
	To be more specific, for a given constant \(a\in\R\)  we are interested in the approximation of \(\psi_{\gamma,T}(u,au)\) as $u\to\IF$.
	Observe that $\psi_{\gamma,T}(u,au)$ can be rewritten in another form, namely
	$$\psi_{\gamma,T}(u,au)=\pk{\exists_{t\in[0,T], s_1,s_2\in[0,t]}:\begin{array}{ccc}B_1(t)-c_1t-\gamma_1(B_1(s_1)-c_1s_1)>u\\B_2(t)-c_2t-\gamma_2(B_2(s_2)-c_2s_2)>au\end{array}}.$$
	
	Using the self-similarity property of the Brownian motion, we can assume without loss of generality that \(T=1\). In this case denote $\psi_{\gamma}(u,au):=\psi_{\gamma,1}(u,au)$. Moreover due to the symmetric view of \(\vk X\) we can assume here and below \(a\leq 1\). 
	
	Here and below bold notation are vector-valued and non-bold designations are scalars. Throughout the whole paper \(t\geq 0\), \(\vk s=(s_1,s_2)\in[0,\IF)^2\), \(\vk\gamma=(\gamma_1,\gamma_2)\in[0,2)^2\). Additionally, hereafter \(\vk B(\vk s)=(B_1(s_1),B_2(s_2))\), \(\vk B(t):=\vk B(t,t)=(B_1(t),B_2(t))\), and $\vk a:=(1,a)$.
	
	All operations for vectors, such as addition, division, comparison and multiplication by scalar constants, are assumed to be componentwise. For example, $\vk x/\vk y:=(x_1/y_1,x_2/y_2)$ for $\vk x, \vk y\in\R^2$.
	
	Thus in terms of the previously discussed vector notation one can rewrite the formula for $\psi_{\gamma}(u,au)$ as
	$$\psi_{\gamma}(u,au)=\psi_{\gamma}(\vk au)=\pk{\exists_{t\in[0,1], s_1,s_2\in[0,t]}:\vk B(t)-\vk ct-\vk\gamma(\vk B(\vk s)-\vk c\vk s)>\vk au}.$$
	
	We formulate next the main theorem.
	\begin{theo} \label{Th1}
		Let $\vk B^*$ be an independent copy of $\vk B$ and let $\phi$ be the density function of the vector $(B_1(1),B_2(1))$. Then for all $1\geq a>0$
		$$\psi_{\gamma}(\vk au)\sim C(a)\frac{1}{u^{2}}\phi(\vk au+\vk c),$$
		where
		$$C(a)=\int_{\R^2}\pk{\exists_{(t,\vk s)\geq \vk 0}: \vk B^*(t)-\vk at+\vk\gamma(\vk B(\vk s)-\vk a\vk s)>\vk x}e^{x_1+ax_2}d\vk x\in(0,\IF).$$
		
		If $a\leq 0$ and $\tilde{\vk a}=(1,0)$, then 
		$$\psi_{\gamma}(\vk au)\sim C(a)\frac{1}{u}\phi(\tilde{\vk a}u+\vk c), \ \quad C(a)=\frac{4}{2-\gamma_1}\sqrt{2\pi}\exp\left(\frac{c_2^2}{2}\right)\Phi^*(-c_2),$$
		with $\Phi^*(x):=1$ for all $x$ if $a<0$ and $\Phi^*(x):=\Phi(x)$ if $a=0$.
	\end{theo}
	
	This paper is organized as follows: in the next section the result of Theorem \ref{Th1} is discussed followed by two illustrating examples. In Section $3$ we present scheme of the proof of Theorem \ref{Th1}. Technical details are postponed to the last section.
	
	\section{Discussions and Examples}
	\begin{remark}\label{Lem2}
		Observe that according to \cite[Lemma 2]{Lemma2} as $u\to\IF$
		\[
		\frac{1}{u^2}\phi(\vk au+\vk c)\sim a\pk{\vk B(1)>\vk au+\vk c}
		\]
		if \(a>0\), and
		$$\sqrt{2\pi}\exp\left(\frac{c_2^2}{2}\right)\Phi^*(-c_2)\frac{1}{u}\phi(\tilde{\vk a}u+\vk c)\sim \pk{\vk B(1)>\vk au+\vk c}$$
		if $a\leq 0$ and $\tilde{\vk a}=(1,0).$ Hence by Theorem \ref{Th1} for \(a>0\) we have that
		\[
		\psi_{\gamma}(\vk au)\sim aC(a)\pk{\vk B(1)>\vk au+\vk c}
		\]
		as \(u\to\IF\), where \(C(a)\) is as in Theorem \ref{Th1}. If \(a\leq 0\), then
		\[
		\psi_{\gamma}(\vk au)\sim \frac{4}{2-\gamma_1}\pk{\vk B(1)>\vk au+\vk c}, \ \ u\to\IF.
		\]
	\end{remark}
	\begin{example}
		Let $\gamma_1=0$. By the Theorem \ref{Th1} we obtain that for $a>0$
		$$\psi_{\gamma}(u,au)\sim C(a)\frac{1}{u^{2}}\phi(\vk au+\vk c)\sim aC(a)\pk{\vk B(1)>\vk au+\vk c},$$
		where
		$$C(a)=\int_{\R^2}\pk{\exists_{(t,s)\geq \vk 0}:\begin{array}{ccc} B^*_1(t)-t>x_1\\B^*_2(t)-at+\gamma_2(B_2(s)-as)>x_2\end{array}}e^{x_1+ax_2}d\vk x.$$
		If $a\leq 0$ and $\tilde{\vk a}=(1,0)$, then 
		$$\psi_{\gamma}(u,au)\sim C(a)\frac{1}{u}\phi(\tilde{\vk a}u+\vk c)\sim2\pk{\vk B(1)>\vk au+\vk c},$$
		where
		$$C(a)=2\sqrt{2\pi}e^{\frac{c_2^2}{2}}\Phi^*(-c_2),$$
		with $\Phi^*(x)$ as in Theorem \ref{Th1}.
		
		In particular, if $\gamma_2=\gamma_1=0$, then for the case $a>0$ the constant $C(a)$ can be simplified:
		$$C(a)=\int_{\R^2}\pk{\exists_{t\geq 0}:\begin{array}{ccc} B_1(t)-t>x_1\\B_2(t)-at>x_2\end{array}}e^{x_1+ax_2}d\vk x.$$
		Notice that the last observation is consistent with the result \cite[Thm 2.1]{mi:18}.
	\end{example}
	\begin{example}
		Let $\gamma_1=\gamma_2=1$. According to Theorem \ref{Th1} for $a>0$
		$$\pk{\exists_{t\in[0,1], s_1,s_2\in[0,t]}:\vk B(t)-\vk ct-\vk B(\vk s)+\vk c\vk s>\vk au}\sim C(a)\frac{1}{u^{2}}\phi(\vk au+\vk c)\sim aC(a)\pk{\vk B(1)>\vk au+\vk c},$$
		where
		$$C(a)=\int_{\R^2}\pk{\exists_{(t,\vk s)\geq \vk 0}: \vk B^*(t)-\vk at+\vk B(\vk s)-\vk a\vk s>\vk x}e^{x_1+ax_2}d\vk x.$$
		If $a\leq 0$ and $\tilde{\vk a}=(1,0)$, then 
		$$\pk{\exists_{t\in[0,1], s_1,s_2\in[0,t]}:\vk B(t)-\vk ct-\vk B(\vk s)+\vk c\vk s>\vk au}\sim C(a)\frac{1}{u}\phi(\tilde{\vk a}u+\vk c)\sim 2\pk{\vk B(1)>\vk au+\vk c},$$
		where
		$$C(a)=2\sqrt{2\pi}e^{\frac{c_2^2}{2}}\Phi^*(-c_2).$$
	\end{example}
	\begin{remark}\label{C(a)}
		If \(a>0\), the constant cannot be calculated explicitly, but one can obtain the following upper bound on \(C(a)\):
		\[
		C(a)\leq \frac{16}{a(2-\gamma_1)(2-\gamma_2)}<\infty.
		\]
	\end{remark}
	\begin{proof}
		Since for the infinite ruin of Brownian motion 
		$$\pk{ \sup_{t\ge 0} (B(t)- ct) > u}= e^{-2c u}, \quad u>0$$
		one can conclude that  $Y_c=\sup_{t\ge 0} (B(t)- ct)$ has exponential distribution with mean $1/(2c)$. 
		Notice that 
		\begin{align*}
			\pk{\exists_{(t,\vk s)\geq \vk 0}: \vk B^*(t)-\vk at+\vk\gamma(\vk B(\vk s)-\vk a\vk s)>\vk x}&\leq\pk{\begin{array}{ccc}\sup_{t_1\ge 0} (B^*(t_1)- t_1)+ \gamma_1\sup_{s_1\ge 0} (B(s_1)- s_1)>x_1\\ \sup_{t_2\ge 0} (B^*(t_2)- at_2)+ \gamma_2\sup_{s_2\ge 0} (B(s_2)- as_2)>x_2\end{array}}\\&=\pk{\begin{array}{ccc}Y^{*}_1+ \gamma_1Y_1>x_1\\ Y^*_a+ \gamma_2Y_a>x_2\end{array}}
			\\&=\pk{Y^*_1+ \gamma_1Y_1>x_1}\pk{Y^*_a+ \gamma_2Y_a>x_2},
		\end{align*}
		where $Y_1,Y_a, Y_1^*, Y_a^*$ are independent exponentially distributed random variables with 
		$$\mathbb{E}\{Y^*_1\}=\mathbb{E}\{Y_1\}=\frac{1}{2} \text{ and } \mathbb{E}\{Y^*_a\}=\mathbb{E}\{Y_a\}=\frac{1}{2a}.$$
		Hence by the independence
		\begin{align*}
			C(a):=\int_{\R^2}\pk{\exists_{(t,\vk s)\geq \vk 0}: \vk B^*(t)-\vk at+\vk\gamma(\vk B(\vk s)-\vk a\vk s)>\vk x}e^{x_1+ax_2}d\vk x&\leq\frac{1}{a}\mathbb{E}\{e^{Y^*_1+\gamma_1Y_1}\}\mathbb{E}\{e^{aY^*_a+\gamma_2aY_a}\}\\&=\frac{1}{a}\mathbb{E}\{e^{Y^*_1}\}\mathbb{E}\{e^{\gamma_1Y_1}\}\mathbb{E}\{e^{aY^*_a}\}\mathbb{E}\{e^{\gamma_2aY_a}\}\\&=\frac{1}{a}\mathbb{E}\{e^{Y^*_1}\}\mathbb{E}\{e^{\gamma_1Y_1}\}\mathbb{E}\{e^{aY^*_a}\}\mathbb{E}\{e^{\gamma_2aY_a}\}\\&=\frac{1}{a}\frac{2}{2-1}\frac{2}{2-\gamma_1}\frac{2a}{2a-a}\frac{2a}{2a-\gamma_2a}\\&=\frac{16}{a(2-\gamma_1)(2-\gamma_2)},
		\end{align*}
		where we used the fact that a moment generating function of an exponential random variable $X$ with $\mathbb{E}\{X\}=\frac{1}{\lambda}$ is equal to $\frac{\lambda}{\lambda-t}$ for $t<\lambda$.
	\end{proof}
	
	\section{Proof of Theorem \ref{Th1}}
	Before we prove Theorem \ref{Th1}, we introduce some auxiliary notation. Denote 
	$$E=\{(x_1,x_2,x_3)\in[0,1]^3:\max(x_2,x_3)\leq x_1\},$$ 
	hence we can rewrite
	$$\psi_{\gamma}(\vk au)=\pk{\exists_{(t,s_1,s_2)\in E}:X_1(t,s_1,s_2)>u,X_2(t,s_1,s_2)>au},$$
	where for \(i=1,2\)
	\begin{equation}\label{X}
		X_i(t,s_1,s_2)=B_i(t)-c_it-\gamma_i(B_i(s_i)-c_is_i), \ \ \ (t,s_1,s_2)\in E.
	\end{equation}
	
	Let $\Lambda, u>0$ and define $$\Delta(u,\Lambda):=\left[1-\frac{\Lambda}{u^2},1\right]\times\left[0,\frac{\Lambda}{u^2}\right]\times\left[0,\frac{\Lambda}{u^2}\right].$$ 
	To simplify notation, we will omit the parameters $u$ and $\Lambda$ and simply write $\Delta$, implying that it depends on $u$ and $\Lambda$. For any fixed \(\Lambda>0\) and large $u$ we have that $\Delta\in E$.
	
	To verify Theorem \ref{Th1}, divide the set $E$ into several parts and address each part separately. For $\Lambda,u>0$ introduce two auxiliary functions
	\begin{align*}
		M(u,\Lambda):=\pk{\exists_{(t,s_1,s_2) \in\Delta}: X_1(t, s_1, s_2)> u, X_2(t, s_1, s_2)> au},
		\\m(u,\Lambda):=\pk{\exists_{(t,s_1,s_2) \in\Delta^c}: X_1(t, s_1, s_2)> u, X_2(t, s_1, s_2)> au},
	\end{align*}
	where for a set $A\subset E$ we denote by $A^c$ the set $E\setminus A,$ with "$E\setminus A$" being the set difference between $E$ and $A$.
	
	It turns out that the main contribution to the asymptotics of $\psi_{\gamma}(\vk au)$ comes from $M(u,\Lambda)$, while $m(u,\Lambda)$ is irrelevant for the asymptotics.

	\BEL
	1. Let $\phi$ be the density funtion of $(B_1(1), B_2(1))$. If $a>0$, then 
	\label{M(u,Lambda)}
	\BQN
	M(u,\Lambda)&\sim &C(a, \Lambda)\frac{1}{u^{2}}\phi(\vk au+\vk c), \quad u\to \IF,
	\EQN
	where
	\BQN
	C(a, \Lambda)=\int_{\R^2}\pk{\exists_{(t,\vk s) \in [0, \Lambda]^3}: \vk B^*(t)-\vk at+\vk\gamma(\vk B(\vk s)-\vk a\vk s)>\vk x}e^{x_1+ax_2}d\vk x\in(0,\IF),
	\EQN
	and $\vk B^*$ is an independent copy of $\vk B.$
	
	2. If $a=0$, then
	\BQN
	M(u,\Lambda)&\sim &C(a, \Lambda)\frac{1}{u}\phi(\vk au+\vk c), \quad u\to \IF,
	\EQN
	where the constant \(C(a, \Lambda)\) is equal to
	\BQN
	\int_{\R^2}\mathbb{I}(x_{2}<0)\pk{\exists_{(t,s) \in [0, \Lambda]^2}: B^*_1(t)-t+\gamma_1(B_1(s)-s)>x_1}e^{x_1-\frac{x_{2}^2-2c_{2}x_{2}}{2}}d\vk x\in(0,\IF).
	\EQN
	If $a<0$, then 
	\BQN
	M(u,\Lambda)&\sim &C(a, \Lambda)\frac{1}{u}\phi(\tilde{\vk a}u+\vk c), \ \quad u\to \IF,
	\EQN
	where
	\BQN
	C(a, \Lambda)=\int_{\R^2}\pk{\exists_{(t,s) \in [0, \Lambda]^2}: B^*_1(t)-t+\gamma_1(B_1(s)-s)>x_1}e^{x_1-\frac{x_{2}^2-2c_{2}x_{2}}{2}}d\vk x\in(0,\IF).
	\EQN
	\EEL
	For $\ve\in(0,1)$ denote by \(A_{\ve}\) the following set
	\[
	A_{\ve}=\{(x_1,x_2,x_3)\in[0,1]^3:\max(x_2,x_3)\leq x_1\leq\ve\}\subset E.
	\]
	Introduce, in addition, the set \(A_0\) which is defined by
	$$A_0(\Lambda,u):=\left(\left[1-\frac{\Lambda\ln^2u}{u^2},1\right]\times\left[0,\frac{\Lambda\ln^2u}{u^2}\right]\times\left[0,\frac{\Lambda\ln^2u}{u^2}\right]\right)\setminus\Delta.$$ Observe that $A_0(\Lambda,u)\subset E$ for large $u$. For simplicity we will write $A_0$. It follows that 
	$$m(u,\Lambda)\leq \pk{\exists_{(t, \vk s)\in A_{\ve}}:\vk X(t,\vk s)>\vk au}+\pk{\exists_{(t, \vk s)\in A_0}:\vk X(t,\vk s)>\vk au}+\pk{\exists_{(t, \vk s)\in A_1}:\vk X(t,\vk s)>\vk au},$$
	where $A_1:=A_1(\Lambda,u):=(\Delta\cup A_0\cup A_{\ve})^c$ and \(\vk X\) is given by \eqref{X}. Thus it is enough to estimate separately 
	$$\pk{\exists_{(t, \vk s)\in A_{\ve}}:\vk X(t,\vk s)>\vk au}, \ \ \ \ \pk{\exists_{(t, \vk s)\in A_0}:\vk X(t,\vk s)>\vk au} \ \ \ \text{ and } \ \ \ \pk{\exists_{(t, \vk s)\in A_1}:\vk X(t,\vk s)>\vk au}.$$
	\BEL
	$i)$ For any $\delta>0$ there exists $\ve=\ve(\delta)>0$ such that
	\label{m(u,Lambda)}
	\BQN\label{A_eps}
	\pk{\exists_{(t,\vk s)\in A_{\ve}}:\vk B(t)-\vk ct-\vk\gamma(\vk B(\vk s)-\vk c\vk s)>\vk au}&=&o\left(\exp\left(-\delta u^2\right)\right), \quad u\to \IF.
	\EQN
	$ii)$ Moreover,
	\BQN\label{A_1}
	\pk{\exists_{(t,\vk s)\in A_1}:\vk B(t)-\vk ct-\vk\gamma(\vk B(\vk s)-\vk c\vk s)>\vk au}=o(\pk{\vk B(1)-\vk c>\vk au}), \quad u\to \IF.
	\EQN
	$iii)$ There exist positive constants $G_1$ and $G_2$ such that
	\BQN\label{A_0}
	\pk{\exists_{(t,\vk s)\in A_0}:\vk B(t)-\vk ct-\vk\gamma(\vk B(\vk s)-\vk c\vk s)>\vk au}\leq G_1e^{-G_2\Lambda}\pk{\vk B(1)-\vk c>\vk au}, \quad u\to \IF.
	\EQN
	\EEL
	Observe that for $\delta>(1+a^2)/2$ we have that \(\exp\left(-\delta u^2\right)=o(\pk{\vk B(1)>\vk au+\vk c})\).
	
	Assume here and below $\delta=1+a^2$. As it was shown above, there exists $\ve=\ve(\delta)$ such that
	$$\pk{\exists_{(t,\vk s)\in A_{\ve}}:\vk B(t)-\vk ct-\vk\gamma(\vk B(\vk s)-\vk c\vk s)>\vk au}=o(\pk{\vk B(1)>\vk au+\vk c}),\ \ u\to\IF.$$
	By the results of Lemma \ref{m(u,Lambda)} 
	$$\lim_{\Lambda\to\IF}\lim_{u\to\IF}\frac{m(u,\Lambda)}{\pk{\vk B(1)>\vk au+\vk c}}=0.$$
	Using Lemma \ref{M(u,Lambda)} and observations discussed in Remark \ref{Lem2} one can therefore conclude that
	$$\lim_{\Lambda\to\IF}\lim_{u\to\IF}\frac{m(u,\Lambda)}{\psi_{\gamma}(\vk au)}=0.$$
	Since $M(u,\Lambda)\leq\psi_{\gamma}(\vk au)\leq M(u,\Lambda)+m(u,\Lambda)$,
	$$\lim_{\Lambda\to\IF}\lim_{u\to\IF}\frac{M(u,\Lambda)}{\psi_{\gamma}(\vk au)}=1.$$
	To complete the proof, we apply the result of Lemma \ref{M(u,Lambda)} together with the monotone convergence theorem to obtain the desired asymptotics of 
	\(\psi_{\gamma}(\vk au)\). For the reader's convenience, we consider the cases \(a>0\) and \(a\leq 0\) separately. If \(a>0\), then by Lemma \ref{M(u,Lambda)} 
	$$M(u,\Lambda)\sim C(a, \Lambda)\frac{1}{u^{2}}\phi(\vk au+\vk c), \quad u\to \IF,$$
	$$C(a, \Lambda)=\int_{\R^2}\pk{\exists_{(t,\vk s)\in [0,\Lambda]^3}: \vk B^*(t)-\vk at+\vk\gamma(\vk B(\vk s)-\vk a\vk s)>\vk x}e^{x_1+ax_2}d\vk x.$$
	By the monotone convergence theorem 
	$$\lim_{\Lambda\to\IF}C(a, \Lambda)=\int_{\R^2}\pk{\exists_{(t,\vk s)\geq 0}: \vk B^*(t)-\vk at+\vk\gamma(\vk B(\vk s)-\vk a\vk s)>\vk x}e^{x_1+ax_2}d\vk x=C(a).$$
	
	By Remark \ref{C(a)} we have that $C(a)<\infty$, hence the statement of Theorem \ref{Th1} follows for \(a>0.\) If $a\leq 0$, then as in the case $a>0$
	$$\lim_{\Lambda\to\IF}C(a, \Lambda)=\int_{\R^2}\pk{\exists_{(t,s) \geq 0}: B^*_1(t)-t+\gamma_1(B_1(s)-s)>x_1}e^{x_1-\frac{x_{2}^2-2c_{2}x_{2}}{2}}d\vk x=:C(a)$$
	if $a<0$ and
	$$\lim_{\Lambda\to\IF}C(a, \Lambda)=\int_{\R^2}\mathbb{I}(x_2<0)\pk{\exists_{(t,s) \geq 0}: B^*_1(t)-t+\gamma_1(B_1(s)-s)>x_1}e^{x_1-\frac{x_{2}^2-2c_{2}x_{2}}{2}}d\vk x=:C(a)$$
	if $a=0$. If $a=0$, then we have that
	\begin{align*}
		&\int_{\R^2}\mathbb{I}(x_2<0)\pk{\exists_{(t,s) \geq 0}: B^*_1(t)-t+\gamma_1(B_1(s)-s)>x_1}e^{x_1-\frac{x_{2}^2-2c_{2}x_{2}}{2}}d\vk x&
		\\&=e^{\frac{c_2^2}{2}}\sqrt{2\pi}\Phi(-c_2)\int_{\R}\pk{\exists_{(t,s) \geq 0}: B^*_1(t)-t+\gamma_1(B_1(s)-s)>x}e^{x}dx.&
	\end{align*}

    Finally, observe that
	\[
	\int_{\R}\pk{\exists_{(t,s) \geq 0}: B^*_1(t)-t+\gamma_1(B_1(s)-s)>x}e^{x}dx=\int_{0}^\IF \pk{ Y^*_1 + \gamma_1 Y_1 > x}e^{x}dx,
	\]
    where \(Y_1^*\) and \(Y_1\) are defined as in the proof of Remark \ref{C(a)}. Since $ Y^*_1 + \gamma_1 Y_1$ is almost surely positive and $Y_1^*, Y_1$ are independent, then
	$$\int_{0}^\IF \pk{ Y^*_1 + \gamma_1 Y_1 > x}e^{x}dx=\mathbb{E}\{e^{Y^*_1 + \gamma_1 Y_1}\}=\mathbb{E}\{e^{Y^*_1}\}\mathbb{E}\{e^{\gamma_1 Y_1}\}=\frac{2}{2-1}\frac{2}{2-\gamma_1}=\frac{4}{2-\gamma_1}.$$
	The case $a<0$ is simpler and the proof is therefore omitted. Hence the proof follows.
	
	\qed
	\section{Proofs of auxiliary Lemmas}
	
	\subsection{Proof of Lemma \ref{M(u,Lambda)}}
	Consider first the case $a>0$. Setting 
	$$\vk u_{x}:=\vk au+\vk c-\frac{\vk x}{u}, t_u:=1-\frac{t}{u^2}, \vk s_{u}=(s_{u,1},s_{u,2}):=\left(\frac{s_1}{u^2},\frac{s_2}{u^2}\right)$$
	by the law of total probability, we obtain
	\begin{align*}
		M(u,\Lambda)&=\pk{\exists_{(t,s_1,s_2) \in [0, \Lambda]^3}: X_1(t_u, \vk s_{u})> u, X_2(t_u, \vk s_{u})> au}
		\\&=\frac{1}{u^2}\int_{\R^2}\pk{\exists_{(t,\vk s) \in [0, \Lambda]^3}: \vk X(t_u, \vk s_{u})> \vk au \,\middle|\, \vk B(1)=\vk u_{x}}\phi(\vk u_{x})d\vk x,
	\end{align*}
	with $\phi$ being the density function of the vector $\vk B(1)$. Since for all $\vk x\in\R^2$ 
	$$\phi(\vk u_{x})=\phi(\vk au+\vk c)\theta_u(\vk x), \ \theta_u(\vk x)=\exp\left(x_1+ax_2+\frac{c_1x_1+c_2x_2}{u}-\frac{x_1^2+x_2^2}{2u^2}\right),$$
	we have
	$$M(u,\Lambda)=\frac{1}{u^2}\phi(\vk au+\vk c)\int_{\R^2}\pk{\exists_{(t,\vk s) \in [0, \Lambda]^3}: \vk X(t_u, \vk s_{u})> \vk au \,\middle|\, \vk B(1)=\vk u_{x}}\theta_u(\vk x)d\vk x.$$
	
	Since $\vk B$ has independent components,
	$$\left(\vk B(t_u) \,\middle|\, \vk B(1)=\vk u_{x}\right)=\vk u_{x}t_u+\vk D(t_u), \ \ \ \ \ \ \ \ \ \ \ \ \ \ \ \ \left(\vk B(\vk s_u) \,\middle|\, \vk B(1)=\vk u_{x}\right)=\vk u_{x}\vk s_u+\vk D(\vk s_u), \,$$
	where \(\vk D\) is defined by
	\begin{equation}\label{D}
		\vk D(\vk v)=\vk B(\vk v)-\vk v\vk B(1), \ \ \forall\vk v\in[0,\infty)^2.
	\end{equation}
	Consequently we have for the conditioning above
	$$\vk u_{x}t_u+\vk D(t_u)-\vk ct_u-\vk\gamma(\vk u_{x}\vk s_{u}+\vk D(\vk s_{u})-\vk c\vk s_{u})>\vk au,$$
	or equivalently
	$$u\vk D(t_u)-u\vk\gamma\vk D(\vk s_{u})-\vk a(t+\vk\gamma\vk s)+\vk f(\vk x, u, t, \vk s)>\vk x,\ \vk f(\vk x, u, t, \vk s)=\frac{t\vk x}{u^2}+\frac{\vk\gamma\vk s\vk x}{u^2},$$
	where hereafter \(t+\vk\gamma\vk s=\vk 1t+\vk\gamma\vk s\), with \(\vk 1=(1, 1).\) 
	One can check that as $u\to\IF$
	$$\vk X_u(\vk x, t, \vk s):=u\vk D(t_u)-u\vk\gamma\vk D(\vk s_{u})-\vk a(t+\vk\gamma\vk s)+\vk f(\vk x, u, t, \vk s)-\vk x$$
	weakly converges to 
	$$(\vk B^*(t)-\vk at)+\vk\gamma(\vk B(\vk s)-\vk a\vk s)-\vk x,$$
	where $\vk B^*$ is an independent copy of $\vk B$.
	Indeed, as $u\to\IF$
	$$\mathbb{E}\left\{u\vk D(t_u)-u\vk\gamma\vk D(\vk s_{u})-\vk a(t+\vk\gamma\vk s)+\vk f(\vk x, u, t, \vk s)-\vk x\right\}\to\mathbb{E}\{\vk B^*(t)-\vk at+\vk\gamma(\vk B(\vk s)-\vk a\vk s)-\vk x\}$$
	for all $\vk x\in\R^2$.
	
	Moreover, if $t_1,t_2,s_{1,1},s_{1,2},s_{2,1},s_{2,2}\in [0,\Lambda]$, then for $i=1,2$ as $u\to\IF$
	$$cov\left(uD_i\left(1-\frac{t_1}{u^2}\right)-u\gamma_iD_i\left(\frac{s_{i,1}}{u^2}\right),uD_i\left(1-\frac{t_2}{u^2}\right)-u\gamma_iD_i\left(\frac{ s_{i,2}}{u^2}\right)\right)\to \min(t_1,t_2)+\gamma_i^2\min(s_{i,1},s_{i,2}).$$
	Observe that since the process $D_i(t)$ does not depend on $\vk B(1)$, one can conclude that
	\begin{align*}
		&u^2cov\left(D_i\left(1-\frac{t_1}{u^2}\right)-\gamma_iD_i\left(\frac{s_{i,1}}{u^2}\right),D_i\left(1-\frac{t_2}{u^2}\right)-\gamma_iD_i\left(\frac{ s_{i,2}}{u^2}\right)\right)&
		\\&=u^2cov\left(D_i\left(1-\frac{t_1}{u^2}\right)-\gamma_iD_i\left(\frac{s_{i,1}}{u^2}\right),B_i\left(1-\frac{t_2}{u^2}\right)-\gamma_iB_i\left(\frac{ s_{i,2}}{u^2}\right)\right).&
	\end{align*}
	Notice that as $u\to\IF$
	$$u^2cov\left(D_i\left(1-\frac{t_1}{u^2}\right), B_i\left(1-\frac{t_2}{u^2}\right)\right)=\min(t_1,t_2)\left(1-\frac{\max(t_1,t_2)}{u^2}\right)\to \min(t_1,t_2),$$
	$$u^2cov\left(B_i\left(\frac{s_{i,1}}{u^2}\right)-\frac{s_{i,1}}{u^2}B_i(1), B_i\left(1-\frac{t_{2}}{u^2}\right)\right)=s_{i,1}-s_{i,1}\left(1-\frac{t_{2}}{u^2}\right)=\frac{s_{i,1}t_{2}}{u^2}\to 0,$$
	$$u^2cov\left(B_i\left(1-\frac{t_1}{u^2}\right)-\left(1-\frac{t_1}{u^2}\right)B_i(1), B_i\left(\frac{s_{i,2}}{u^2}\right)\right)=s_{i,2}-s_{i,2}\left(1-\frac{t_{1}}{u^2}\right)=\frac{s_{i,2}t_{1}}{u^2}\to 0$$
	and
	$$cov\left(u\gamma_iD_i\left(\frac{s_{i,1}}{u^2}\right),u\gamma_iB_i\left(\frac{s_{i,2}}{u^2}\right)\right)=u^2\gamma_i^2\left(\frac{\min(s_{i,1},s_{i,2})}{u^2}\right)-u^2\gamma_i^2\frac{s_{i,1}s_{i,2}}{u^4}\to \gamma_i^2\min(s_{i,1},s_{i,2}).$$
	Finally, we check the tightness. For $i=1,2$ we prove
	$$\mathbb{E}\{(X_{u,i}(\x, t_1, \vk s_1)-X_{u,i}(\x, t_2, \vk s_2))^2\}\leq C(|t_1-t_2|+|s_{1,1}-s_{1,2}|+|s_{2,1}-s_{2,2}|)$$
	for some constant $C$.
	By the triangle inequality it is sufficient to check that
	\begin{align*}
		1. \ &  \mathbb{E}\left\{\left(uD_i\left(1-\frac{t_1}{u^2}\right)-uD_i\left(1-\frac{t_2}{u^2}\right)\right)^2\right\}\leq C_1|t_1-t_2|,
		\\ 2. \ &  \mathbb{E}\left\{\left(uD_i\left(\frac{s_{i,1}}{u^2}\right)-uD_i\left(\frac{s_{i,2}}{u^2}\right)\right)^2\right\}\leq C_2|s_{i,1}-s_{i,2}|,
		\\ 3. \ & \left(a_i\left(t_1+\gamma_is_{i,1}\right)-a_i\left(t_2+\gamma_is_{i,2}\right)\right)^2\leq C_3(|t_1-t_2|+|s_{i,1}-s_{i,2}|),
		\\ 4. \ & \left(f\left(x_i, u, t_1, s_{i,1}\right)-f\left(x_i, u, t_2, s_{i,2}\right)\right)^2\leq C_4(x_i)(|t_1-t_2|+|s_{i,1}-s_{i,2}|)
	\end{align*}
	
	for positive constants $C_1, C_2, C_3, C_4(x_i)$ and $t_1, t_2, s_{1,1},s_{2,1}s_{1,2},s_{2,2}\in [0,\Lambda]$.
	
	The inequalities $3$ and $4$ are obvious, and the first two follow since
	for $t_1,t_2\in [0,1]$ we have 
	\begin{align*}
		\mathbb{E}\{(D_i(t_1)-D_i(t_2))^2\}&=cov(D_i(t_1)-D_i(t_2),B_i(t_1)-B_i(t_2))\\
		&=cov(B_i(t_1)-B_i(t_2)-t_1B_i(1)+t_2B_i(1),B_i(t_1)-B_i(t_2))\\
		&=t_1-\min(t_1,t_2)-t_1^2+t_1t_2-\min(t_1,t_2)+t_2+t_1t_2-t_2^2\\
		&=|t_1-t_2|(1-|t_1-t_2|)<|t_1-t_2|.
	\end{align*}
	Consequently the claimed weak convergence holds.
	
	Since the $\sup$ functional is continuous on the space $C([0,\Lambda]^3)\times C([0,\Lambda]^3)$, by the continuous mapping theorem
	\begin{align*}
		\lim_{u\to\IF}\pk{\sup_{(t,\vk s) \in [0, \Lambda]^3}\vk X_u(\vk x, t,\vk s)>\vk 0}
		&=\pk{\sup_{(t,\vk s) \in [0, \Lambda]^3}(\vk B^*(t)-\vk at+\vk\gamma(\vk B(\vk s)-\vk a\vk s)-\vk x)>\vk 0}
	\end{align*}
	for almost all $\vk x$. Moreover, we have \(\lim_{u\to\IF}\theta_u(\vk x)=\exp(x_1+ax_2).\)

	Our goal is to use a dominated convergence theorem and for that one need to estimate $$\pk{\exists_{(t,\vk s)\in[0,\Lambda]^3}: \vk X_u(\vk x, t, \vk s)>\vk 0}.$$
	Observe that for large enough $u$
	\begin{align*}
		\pk{\exists_{(t,\vk s)\in[0,\Lambda]^3}:\vk X_u(\vk x, t, \vk s)>\vk 0}&\leq\pk{\exists_{(t,\vk s)\in[0,\Lambda]^3}:u\vk D\left(t_u\right)-u\vk\gamma\vk D\left(\vk s_{u}\right)>G_1\vk x}\\
		&\leq\pk{\exists_{(t,\vk s)\in[0,\Lambda]^3}: u\sum_{i=1}^{2}\left(D_i\left(t_u\right)-\gamma_iD_i\left(\vk s_{u,i}\right)\right)>G_1\sum_{i=1}^{2}x_i}
	\end{align*}
	for $G_1>0$.
	Notice that since $B_1$ and $B_2$ are independent
	$$Var\left(uD_1\left(t_u\right)-u\gamma_1D_1\left(s_{u,1}\right)+uD_2\left(t_u\right)-u\gamma_2D_2\left(\vk s_{u,2}\right)\right)$$ is equal to $$Var\left(uD_1\left(t_u\right)-u\gamma_1D_1\left(s_{u,1}\right)\right)+Var\left(uD_2\left(t_u\right)-u\gamma_2D_2\left(s_{u,2}\right)\right).$$
	Observe that for $i=1,2$ we have that
	\begin{align*}
		Var\left(uD_i\left(t_u\right)-u\gamma_iD_i\left(s_{u,i}\right)\right)&=cov\left(uD_i\left(t_u\right)-u\gamma_iD_i\left(\vk s_{u,i}\right),uD_i\left(t_u\right)-u\gamma_iD_i\left(\vk s_{u,i}\right)\right)
		\\&=cov\left(uD_i\left(t_u\right)-u\gamma_iD_i\left(s_{u,i}\right),uB_i\left(t_u\right)-u\gamma_iB_i\left(s_{u,i}\right)\right)\\
		&=t\left(1-\frac{t}{u^2}\right)-\frac{2s_i\gamma_it}{u^2}+s_i\gamma_i\left(1-\frac{s_i}{u^2}\right)\in \left(\frac{t+s_i\gamma_i}{2}, t+s_i\gamma_i\right).
	\end{align*}
	
	Using Borel's inequality (see \cite[Lemma 4.5]{Vector}), we can apply the dominated convergence theorem to $\vk X_u(\vk x, t,\vk s)$ and conclude that
	$$M(u,\Lambda)\sim\frac{1}{u^2}\phi(\vk au+\vk c) \int_{\R^2}\pk{\exists_{(t,\vk s) \in [0, \Lambda]^3}: \vk B^*(t)-\vk at+\vk\gamma(\vk B(\vk s)-\vk a\vk s)>\vk x}e^{x_1+ax_2}d\vk x.$$
	Consider now the case $a\leq 0$ and denote $\vk u_{x}=\tilde{\vk a}u+\vk c-\vk x/\tilde{\vk u}$, where $\tilde{\vk u}=(u,1)$ and $\tilde{\vk a}=(1,0)$. Then similarly to the case \(a>0\), by the law of total probability, we obtain that
	$$M(u,\Lambda)=\frac{1}{u}\int_{\R^2}\pk{\exists_{(t,\vk s) \in [0, \Lambda]^3}: \vk X\left(t_u, \vk s_{u}\right)> \vk au \,\middle|\, \vk B(1)=\vk u_{x}}\phi(\vk u_{x})d\vk x.$$
	Notice that in this case 
	$$\phi(\vk u_{x})=\phi(\tilde{\vk a}u+\vk c)\theta_u(\vk x),\ \ \theta_u(\vk x)=\exp\left(x_1-\frac{x_2^2-2c_2x_2}{2}+\frac{c_1x_1}{u}-\frac{x_1^2}{2u^2}\right).$$
	Therefore
	$$M(u,\Lambda)=\frac{1}{u}\phi(\tilde{\vk a}u+\vk c)\int_{\R^2}\pk{\exists_{(t,\vk s) \in [0, \Lambda]^3}: \vk X\left(t_u, \vk s_{u}\right)> \vk au \,\middle|\, \vk B(1)=\vk u_{x}}\theta_u(\vk x)d\vk x.$$
	The condition \(\left(\vk X\left(t_u, \vk s_{u}\right)> \vk au \,\middle|\, \vk B(1)=\vk u_{x}\right)\)
	can be rewritten as
	$$t_u\vk u_{x}+\vk D\left(t_u\right)-\vk ct_u-\vk\gamma\left(\vk s_{u}\vk u_{x}+\vk D\left(\vk s_{u}\right)-\vk c\vk s_{u}\right)>\vk au,$$
	where \(\vk D\) is defined by \eqref{D}. That is equivalent to
	$$u\vk D\left(t_u\right)-u\vk\gamma\vk D\left(\vk s_{u}\right)-\tilde{\vk a}(t+\vk\gamma\vk s)+\vk f(\vk x, u, t, \vk s)>\frac{u\vk x}{\tilde{\vk u}}+(\vk a-\tilde{\vk a})u^2,\ \vk f(\vk x, u, t, \vk s)=\frac{t\vk x}{\tilde{\vk u}u}+\frac{\vk\gamma\vk s\vk x}{\tilde{\vk u}u}.$$
	As in the case $a>0$
	$$\vk X_u(\vk x, t, \vk s):=u\vk D\left(t_u\right)-u\vk\gamma\vk D\left(\vk s_{u}\right)-\tilde{\vk a}(t+\vk\gamma\vk s)+\vk f(\vk x, u, t, \vk s)-\frac{u\vk x}{\tilde{\vk u}}-(\vk a-\tilde{\vk a})u^2$$
	weakly converges as $u\to\IF$ to
	\begin{align*}
		1. \ & \ ((B^*_1(t)-t)+\gamma_1(B_1(s_1)-s_1)-x_1, \IF), & \text{ if } a<0 & \text{ or } a=0, x_2<0,
		\\2. \ & \ ((B^*_1(t)-t)+\gamma_1(B_1(s_1)-s_1)-x_1, B^*_2(t)+\gamma_2B_2(s_2)), & \text{ if } a=0 & \text{ and } x_2=0,
		\\3. \ & \ ((B^*_1(t)-t)+\gamma_1(B_1(s_1)-s_1)-x_1, -\IF), & \text{ if } a=0 & \text{ and } x_2>0,
	\end{align*}
	
	where $\vk B^*$ is an independent copy of $\vk B.$ By the continuous mapping theorem 
	\[\pk{\exists_{(t,\vk s) \in [0, \Lambda]^3}: \vk X_u(\vk x, t, \vk s)>\vk au}\] converges as $u\to\IF$ to 
	\begin{align*}
		1. \ & \ \pk{\exists_{(t,\vk s) \in [0, \Lambda]^3}: B^*_1(t)-t+\gamma_1(B_1(s_1)-s_1)>x_1}, & \text{ if } & a<0 \text{ or } a=0, x_2<0,
		\\2. \ & \ \pk{\exists_{(t,\vk s) \in [0, \Lambda]^3}: B^*_1(t)-t+\gamma_1(B_1(s_1)-s_1)>x_1, B^*_2(t)+\gamma_2B_2(s_2)>0}, & \text{ if } & a=0, \ x_2=0,
		\\3. \ & \ 0, & \text{ if } & a=0, \ x_2>0.
	\end{align*}
	
	As in the case $a>0$ one can obtain that there exist positive constants $G_4,G_5$ such that
	$$\pk{\exists_{(t,\vk s)\in[0,\Lambda]^3}:\vk X_u(\vk x, t, \vk s)>\vk 0}\leq G_4e^{-G_5x_1^2}.$$
	Observe that there exist positive $H_0, H_1, b_1(x_1),b_2(x_2)$ such that
	$$\pk{\exists_{(t,\vk s)\in[0,\Lambda]^3}:\vk X_u(\vk x, t, \vk s)>\vk 0}\theta_u(\vk x)\leq e^{b_1(x_1)x_1+b_2(x_2)}h(x_1, x_2),$$
	where $$h(x_1,x_2)=\begin{cases} 0, x_2>0;\\H_0e^{-H_1x_1^2}, x_2\leq 0,\end{cases}, \ b_1(x_1)=\begin{cases} 2, x_1\geq 0;\\\frac{1}{2}, x_1<0\end{cases}, \ b_2(x_2)=-\frac{x_2^2-2c_2x_2}{2}.$$
	
	The function $e^{b_1(x_1)x_1+b_2(x_2)}h_1(x_1)h_2(x_2)$ is integrable. Hence one can apply dominated convergence theorem and obtain that 
	$$M(u,\Lambda)\sim\frac{1}{u}\phi(\tilde{\vk a}u+\vk c) \int_{\R^2}\mathbb{I}(x_2<0)\pk{\exists_{(t,s) \in [0, \Lambda]^2}: B^*_1(t)-t+\gamma_1(B_1(s)-s)>x_1}e^{x_1-\frac{x_2^2-2c_2x_2}{2}}d\vk x$$
	for the case $a=0$
	and 
	$$M(u,\Lambda)\sim\frac{1}{u}\phi(\tilde{\vk a}u+\vk c) \int_{\R^2}\pk{\exists_{(t,s) \in [0, \Lambda]^2}: B^*_1(t)-t+\gamma_1(B_1(s)-s)>x_1}e^{x_1-\frac{x_2^2-2c_2x_2}{2}}d\vk x,$$
	if $a<0$.
	
	Finally, notice that according to \cite[Lemma 2]{Lemma2} for large $u$ we have
	$$\frac{1}{2au^2}\phi(\vk au+\vk c)\leq\pk{\vk B(1)-\vk c>\vk au}\leq M(u,\Lambda),\ a>0,$$
	$$\frac{\sqrt{2\pi}\exp\left(c_2^2/2\right)\Phi^*(-c_2)}{2u}\phi(\tilde{\vk a}u+\vk c)\leq\pk{\vk B(1)-\vk c>\vk au}\leq M(u,\Lambda),\ a\leq 0.$$
	Thus $0<C(a, \Lambda)\leq C(a)<\IF$
	for $\Lambda>0$, where the upper bound for $C(a, \Lambda)$ was already checked in the Section $2$. 
	
	\qed
	
	Next our goal is to prove Lemma \ref{m(u,Lambda)} and in the following three subsections we prove asymptotic estimates (\ref{A_eps}), (\ref{A_1}) and (\ref{A_0}).
	
	\subsection{Proof of (\ref{A_eps})}
	Recall that 
	$$A_{\ve}=\{(x_1,x_2,x_3)\in[0,1]^3:\max(x_2,x_3)\leq x_1\leq\ve\}.$$ Let also
	$$A_{\ve,1}=\{(x_1,x_2)\in[0,1]^2:x_2\leq x_1\leq\ve\}.$$
	Observe that 
	$$\pk{\exists_{(t,\vk s)\in A_{\ve}}:\vk B(t)-\vk ct-\vk\gamma(\vk B(\vk s)-\vk c\vk s)>\vk au}\leq\pk{\exists_{(t,s)\in A_{\ve,1}}:B_1(t)-c_1t-\gamma_1(B_1(s)-c_1s)>u}.$$
	Thus it is sufficient to focus only on the first coordinate of the process and to make designations simpler further in the proof we omit everywhere the subscript $1$. Denote 
	$$X_u(t,s)=\frac{B(t)-\gamma B(s)}{u+ct-\gamma cs}$$
	with variance $\sigma^2_u(t,s)$. Observe that $X_u(t,s)$ is a centered Gaussian process with continuous trajectories. Thus by Borel's inequality (Lemma $4.5$ from \cite{Vector}) for any $\ve>0$ there exists $\mu_{\ve, u}$ such that for $v>\mu_{\ve, u}$ 
	$$\pk{\exists_{(t,s)\in A_{\ve,1}}:X_u(t,s)>bv}\leq e^{-\frac{(v-\mu_{\ve, u})^2}{2\sigma_{u}^2}},$$
	where
	$$\sigma_{u}^2:=\sup_{(t,\vk s)\in A_{\ve,1}}\sigma^2_u(t,s)<\IF \text{ and } \mu_{\ve, u}=\mathbb{E}\left\{\sup_{(t,s)\in A_{\ve,1}}X_u(t,s)\right\}.$$ 
	Observe that for large enough $u$
	$$\sigma_{u}^2=\sup_{(t,s)\in A_{\ve,1}}\left(\frac{t-2\gamma s+\gamma^2s^2}{(u+ct-\gamma cs)^2}\right)\leq\sup_{(t,s)\in A_{\ve,1}}\left(\frac{2(t-2\gamma s+\gamma^2s^2)}{u^2}\right)\leq\frac{2\ve}{u^2}<\IF.$$
	Now let us estimate $\mu_{\ve, u}$ (for large $u$):
	\begin{align*}
		&\mu_{\ve, u}=\mathbb{E}\left\{\sup_{(t,s)\in A_{\ve,1}}X_u(t,s)\right\}\leq \mathbb{E}\left\{\sup_{(t,s)\in A_{\ve,1}}\frac{B(t)-\gamma B(s)}{(u+ct-\gamma cs)}\right\}&
		\\&\leq\mathbb{E}\left\{\sup_{(t,s)\in A_{\ve,1}}\frac{B(t)}{(u+ct-\gamma cs)}\right\}+\gamma\mathbb{E}\left\{\sup_{(t,s)\in A_{\ve,1}}\frac{ B(s)}{(u+ct-\gamma cs)}\right\}&
		\\&\leq \mathbb{E}\left\{\sup_{(t,s)\in A_{\ve,1}}\frac{2B(t)}{u}\right\}+\gamma\mathbb{E}\left\{\sup_{(t,s)\in A_{\ve,1}}\frac{2B(s)}{u}\right\}\leq\frac{4}{u}\mathbb{E}\left\{\sup_{0\leq t\leq \ve}B(t)\right\}.&
	\end{align*}
	
	Since $\sup_{0\leq t\leq\ve}B(t)$ has the same distribution as $|B(\ve)|$,
	then we obtain
	\begin{align*}
		\mathbb{E}\left\{\sup_{0\leq t\leq\ve}B(t)\right\}=\mathbb{E}\left\{|B(\ve)|\right\}\leq\sqrt{\mathbb{E}\{B^2(\ve)\}}=\sqrt{\ve},
	\end{align*}
	where we used the Hölder's inequality.
	
	Therefore there exists a constant $C>0$ such that $\mu_{\ve, u}\leq C\frac{\sqrt{\ve}}{u}\leq 1,$
	for large enough $u$ and $\mu_{\ve, u}\to 0$ as $u\to\IF$. Hence for large enough $u$
	\begin{align*}
		\pk{\exists_{(t,s)\in A_{\ve,1}}:B(t)-ct-\gamma(B(s)-cs)>u}=\pk{\exists_{(t,s)\in A_{\ve,1}}:\vk X_u(t,\vk s)>1}\leq e^{-\frac{(1-\mu_{\ve, u})^2u^2}{2\ve}}\leq e^{-\frac{u^2}{f_u(\ve)}},
	\end{align*}
	where $f_u(\ve)$ is a function, which converges to $0$ uniformly on $u$, as $\ve\to 0$.
	
	Therefore there exists $u_0>0$ and $\ve=\ve(\delta)>0$ such that for any $u\geq u_0$
	$$\pk{\exists_{(t,\vk s)\in A_{\ve}}:\vk B(t)-\vk ct-\vk\gamma(\vk B(\vk s)-\vk c\vk s)>\vk au}=o(e^{-\delta u^2}).$$
	
	\subsection{Proof of (\ref{A_1})}
    Consider first the case $\gamma_1,\gamma_2>0$. Then notice that \(A_1=\cup_{j=1}^{3} A_{1,j}\), where
    \[
    A_{1,1}:=\left[0, 1-\frac{\Lambda\ln^2u}{u^2}\right]\times\left[0, 1\right]^2\cap A_{\ve}^c,
    \]
    \[
    A_{1,2}:=\left[0,1\right]\times\left[\frac{\Lambda\ln^2u}{u^2},1\right]\times\left[0, 1\right]\cap A_{\ve}^c,
    \]
    \[
    A_{1,3}:=\left[0,1\right]\times\left[0, 1\right]\times\left[\frac{\Lambda\ln^2u}{u^2},1\right]\cap A_{\ve}^c.
    \]
    Hence 
    $$\pk{\exists_{(t,\vk s)\in A_1}:\vk B(t)-\vk ct-\vk\gamma(\vk B(\vk s)-\vk c\vk s)>\vk au}\leq\sum_{j=1}^{3}\pk{\exists_{(t,\vk s)\in A_{1,j}}:\vk B(t)-\vk ct-\vk\gamma(\vk B(\vk s)-\vk c\vk s)>\vk au}.$$
    Thus it is sufficient to prove that
    $$\pk{\exists_{(t,\vk s)\in A_{1,j}}:\vk B(t)-\vk ct-\vk\gamma(\vk B(\vk s)-\vk c\vk s)>\vk au}=o(\pk{\vk B(1)-\vk c>\vk au}), \quad u\to \IF$$
    for $j=1,2,3$. Assume first $a\neq 0$. Then for $u$ large enough
    \begin{align*}
    	\pk{\exists_{(t,\vk s)\in A_{1,j}}:\vk B(t)-\vk ct-\vk\gamma(\vk B(\vk s)-\vk c\vk s)>\vk au}=\pk{\exists_{(t,\vk s)\in A_{1,j}}:\vk X_u(t,\vk s)>\vk au}, 
    \end{align*}
    $$\vk X_u(t,\vk s):=\frac{\vk B(t)-\vk\gamma\vk B(\vk s)}{\vk 1+\frac{\vk c}{\vk au}t-\vk\gamma\frac{\vk c}{\vk au}\vk s}.$$
	Our goal is to apply Piterbarg's inequality (Lemma $4.5$ from \cite{Vector}). To this end, one can see that for $(t,s_1,s_2),(r,q_1,q_2)\in A_1$
	$$\sum_{i=1}^{2}\mathbb{E}\left\{\left(\frac{B_i(t)-\gamma_iB_i(s_i)}{1+\frac{c_i}{ a_iu}t-\gamma_i\frac{c_i}{a_iu}s_i}-\frac{B_i(r)-\gamma_iB_i(q_i)}{1+\frac{c_i}{ a_iu}r-\gamma_i\frac{c_i}{a_iu}q_i}\right)^2\right\}\leq C(|t-r|+|s_1- q_1|+|s_2- q_2|)$$
	with $C>0$ independent of $u$ and $(a_1,a_2):=(1,a)$. Indeed,
	\begin{align*}
		\sum_{i=1}^{2}\mathbb{E}\left\{\left(\frac{B_i(t)-\gamma_iB_i(s_i)}{1+\frac{c_i}{ a_iu}t-\gamma_i\frac{c_i}{a_iu}s_i}-\frac{B_i(r)-\gamma_iB_i(q_i)}{1+\frac{c_i}{ a_iu}r-\gamma_i\frac{c_i}{a_iu}q_i}\right)^2\right\}\\\leq2\sum_{i=1}^{2}\mathbb{E}\left\{\left(\frac{B_i(t)}{1+\frac{c_i}{ a_iu}t-\gamma_i\frac{c_i}{a_iu}s_i}-\frac{B_i(r)}{1+\frac{c_i}{ a_iu}r-\gamma_i\frac{c_i}{a_iu}q_i}\right)^2\right\}\\+2\sum_{i=1}^{2}\mathbb{E}\left\{\left(\frac{\gamma_iB_i(s_i)}{1+\frac{c_i}{ a_iu}t-\gamma_i\frac{c_i}{a_iu}s_i}-\frac{\gamma_iB_i(q_i)}{1+\frac{c_i}{ a_iu}r-\gamma_i\frac{c_i}{a_iu}q_i}\right)^2\right\}\\\leq2\sum_{i=1}^{2}\mathbb{E}\left\{\left(\frac{B_i(t)}{1+\frac{c_i}{ a_iu}t-\gamma_i\frac{c_i}{a_iu}s_i}-\frac{B_i(r)}{1+\frac{c_i}{ a_iu}r-\gamma_i\frac{c_i}{a_iu}q_i}\right)^2\right\}\\+2\sum_{i=1}^{2}\mathbb{E}\left\{\left(\frac{B_i(s_i)}{1+\frac{c_i}{ a_iu}t-\gamma_i\frac{c_i}{a_iu}s_i}-\frac{B_i(q_i)}{1+\frac{c_i}{ a_iu}r-\gamma_i\frac{c_i}{a_iu}q_i}\right)^2\right\}.
	\end{align*}
	We check that
	$$\mathbb{E}\left\{\left(\frac{B_i(t)}{1+\frac{c_i}{ a_iu}t-\gamma_i\frac{c_i}{a_iu}s_i}-\frac{B_i(r)}{1+\frac{c_i}{ a_iu}r-\gamma_i\frac{c_i}{a_iu}q_i}\right)^2\right\}\leq C(|t-r|+|s_1- q_1|+|s_2- q_2|)$$
	for $i=1,2$ and some $C>0$.
	\begin{align*}
		&\mathbb{E}\left\{\left(\frac{B_i(t)}{1+\frac{c_i}{ a_iu}t-\gamma_i\frac{c_i}{a_iu}s_i}-\frac{B_i(r)}{1+\frac{c_i}{ a_iu}r-\gamma_i\frac{c_i}{a_iu}q_i}\right)^2\right\}&
		\\\leq&2\mathbb{E}\left\{\left(B_i(t)\left(1+\frac{c_i}{ a_iu}r-\gamma_i\frac{c_i}{a_iu}q_i\right)-B_i(r)\left(1+\frac{c_i}{ a_iu}t-\gamma_i\frac{c_i}{a_iu}s_i\right)\right)^2\right\}&
		\\&=t\left(1+\frac{c_i}{ a_iu}r-\gamma_i\frac{c_i}{a_iu}q_i\right)^2+r\left(1+\frac{c_i}{ a_iu}t-\gamma_i\frac{c_i}{a_iu}s_i\right)^2&
		\\&-2\min(r,t)\left(1+\frac{c_i}{ a_iu}r-\gamma_i\frac{c_i}{a_iu}q_i\right)\left(1+\frac{c_i}{ a_iu}t-\gamma_i\frac{c_i}{a_iu}s_i\right).&
	\end{align*}
	Without loss of generality assume that $t\geq r$. Then since $t,s_i,r,q_i\in[0,1]$
	\begin{align*}
		&\mathbb{E}\left\{\left(\frac{B_i(t)}{1+\frac{c_i}{ a_iu}t-\gamma_i\frac{c_i}{a_iu}s_i}-\frac{B_i(r)}{1+\frac{c_i}{ a_iu}r-\gamma_i\frac{c_i}{a_iu}q_i}\right)^2\right\}&
		\\&\leq(t-r)\left(1+\frac{c_i}{ a_iu}r-\gamma_i\frac{c_i}{a_iu}q_i\right)^2+r\left(\left(1+\frac{c_i}{ a_iu}r-\gamma_i\frac{c_i}{a_iu}q_i\right)-\left(1+\frac{c_i}{ a_iu}t-\gamma_i\frac{c_i}{a_iu}s_i\right)\right)^2&
		\\&\leq 2|t-r|+2r(t-r)^2+2r(s_i-q_i)^2\leq 3(|t-r|+|s_i-q_i|)&
	\end{align*}
	for large $u$. Using similar arguments, one can check that for sufficiently large $u$
	$$\mathbb{E}\left\{\left(\frac{B_i(s_i)}{1+\frac{c_i}{ a_iu}t-\gamma_i\frac{c_i}{a_iu}s_i}-\frac{B_i(q_i)}{1+\frac{c_i}{ a_iu}r-\gamma_i\frac{c_i}{a_iu}q_i}\right)^2\right\}\leq 3\gamma_i^2(|t-r|+|s_i-q_i|).$$
	Now calculate the covariance matrix of our process
	$$\Sigma_u(t,s)=\begin{pmatrix}
		\frac{t-2\gamma_1s_1+\gamma_1^2s_1}{(1+\frac{c_1}{u}t-\gamma_1\frac{c_1}{u}s_1)^2} & 0 \\
		0 & \frac{t-2\gamma_2s_2+\gamma_2^2s_2}{(1+\frac{c_2}{au}t-\gamma_2\frac{c_2}{au}s_2)^2}
	\end{pmatrix}.$$
	Then 
	$$\Sigma_u^{-1}(t,s)=\begin{pmatrix}
		\frac{(1+\frac{c_1}{u}t-\gamma_1\frac{c_1}{u}s_1)^2}{t-2\gamma_1s_1+\gamma_1^2s_1} & 0 \\
		0 & \frac{(1+\frac{c_2}{au}t-\gamma_2\frac{c_2}{au}s_2)^2}{t-2\gamma_2s_2+\gamma_2^2s_2}
	\end{pmatrix}.$$
	
	Observe that since for \(\gamma_1, \gamma_2\in (0,2)\) functions 
	$$f_i(x_1,x_2)=\frac{(1+\frac{c_i}{a_iu}x_1-\gamma_i\frac{c_i}{a_iu}x_2)^2}{x_1-2\gamma_ix_2+\gamma_i^2x_2}$$
	are smooth on the compact set 
	$$[0,1]^2\cap\{0\leq x_2\leq x_1, x_1\geq\ve\},$$ 
	they are Lipschitz continuous with respect to standard distance in $\R^2$.
	
	Since $\Sigma^{-1}_u(t,\vk s)$ is a diagonal matrix then
	$$||\Sigma^{-1}_u(t,\vk s)-\Sigma^{-1}_u(r,\vk q)||<C(|t-r|+|s_1-q_1|+|s_2-q_2|),$$
	where $C$ can be chosen independent of $u$ and $||\cdot||$ is a standard $l_2$ matrix norm.
	
	Then by \cite[Lemma 4.5]{Vector}
	$$\pk{\exists_{(t,\vk s)\in A_{1,j}}:\vk B(t)-\vk ct-\vk\gamma(\vk B(\vk s)-\vk c\vk s)>\vk au}\leq Cu^3e^{-\frac{u^2}{2\sigma_{a,j}^2}},$$
	where 
	$$\sigma_{a,j}^2:=\max_{(t,\vk s)\in A_{1,j}}\sigma_{a}^2(t,\vk s) \ , \ \sigma_{a}^2(t,\vk s)=\frac{1}{\min_{\vk x\geq\vk a}\vk x^T\Sigma_u^{-1}(t,\vk s)\vk x}$$
	and $C$ is some positive constant.
	
	So, our next aim is to explore the function
	$$f_u(t,\vk s):=\vk x^T\Sigma_u^{-1}(t,\vk s)\vk x=\frac{(1+\frac{c_1}{u}t-\gamma_1\frac{c_1}{u}s_1)^2}{t-2\gamma_1s_1+\gamma_1^2s_1}x_1^2+\frac{(1+\frac{c_2}{au}t-\gamma_2\frac{c_2}{au}s_2)^2}{t-2\gamma_2s_2+\gamma_2^2s_2}x_2^2.$$
	Notice that for large enough $u$ the function $f_u(t,\vk s)$ is increasing as a function of $s_1$ and $s_2$ (when $t$ is any fixed parameter such that $s_i\leq t$ for $i=1,2$) and $\vk x\in\R^2$ is any fixed vector.
	
	Indeed, calculate $\frac{\partial f_u(t,\vk s)}{\partial s_1}$ and $\frac{\partial f_u(t,\vk s)}{\partial s_2}$ and check that they are positive for large $u$
	\begin{align*}
		&\frac{\partial f_u(t,\vk s)}{\partial s_i}=\frac{2(1+\frac{c_i}{a_iu}t-\gamma_i\frac{c_i}{a_iu}s_i)(-\gamma_i\frac{c_i}{a_iu})}{t-2\gamma_is_i+\gamma_i^2s_i}-\frac{(1+\frac{c_i}{a_iu}t-\gamma_i\frac{c_i}{a_iu}s_i)^2(-2\gamma_i+\gamma_i^2)}{(t-2\gamma_is_i+\gamma_i^2s_i)^2}\to\frac{2\gamma_i-\gamma_i^2}{(t-2\gamma_is_i+\gamma_i^2s_i)^2}> 0&
	\end{align*}
	for $i=1,2$ and $(1,a)=(a_1,a_2)$.
	
	Moreover, $f_u(t,\vk 0)$ is decreasing for large $u$ as a function of $t$:
	$$\frac{\partial f_u(t,\vk 0)}{\partial t}=\frac{2\frac{c_i}{a_iu}(1+\frac{c_i}{a_iu}t)t-\left(1+\frac{c_i}{a_iu}t\right)}{t^2}\to-\frac{1}{t^2}<0$$
	for $i=1,2$. Thus for large $u$ the maximal point of $$\sigma_{a}^2(t,\vk s)=\frac{1}{\min_{\vk x\geq\vk a}\vk x^T\Sigma_u^{-1}(t,\vk s)\vk x}=\frac{1}{\tilde{\vk a}^T\Sigma_u^{-1}(t,\vk s)\tilde{\vk a}},$$
	is $(1-\frac{\ln^2u\Lambda}{u^2}, 0, 0)$ for $j=1$, $(1, \frac{\ln^2u\Lambda}{u^2}, 0)$ for $j=2$ and $(1,0,\frac{\ln^2u\Lambda}{u^2})$ for $j=3$.
	
	Here $\tilde{\vk a}=(1,a_+)$, where $a_+:=\max(a,0)$. Thus in our case the generalised variance
	$$\sigma_{a,1}^2=\frac{1}{\frac{(1+\frac{c_1}{u}(1-\frac{\ln^2u\Lambda}{u^2}))^2}{1-\frac{\ln^2u\Lambda}{u^2}}+\frac{(1+\frac{c_2}{au}(1-\frac{\ln^2u\Lambda}{u^2}))^2}{1-\frac{\ln^2u\Lambda}{u^2}}a_+^2}=\frac{1-\frac{\ln^2u\Lambda}{u^2}}{(1+\frac{c_1}{u}(1-\frac{\ln^2u\Lambda}{u^2}))^2+(1+\frac{c_2}{au}(1-\frac{\ln^2u\Lambda}{u^2}))^2a_+^2}$$ for $j=1$,
	$$\sigma_{a,2}^2=\frac{1}{\frac{(1+\frac{c_1}{u}-\gamma_1\frac{c_1}{u}\frac{\ln^2u\Lambda}{u^2})^2}{1-2\gamma_1\frac{\ln^2u\Lambda}{u^2}+\gamma_1^2\frac{\ln^2u\Lambda}{u^2}}+(1+\frac{c_2}{au})^2a_+^2}$$ for $j=2$,
	$$\sigma_{a,3}^2=\frac{1}{(1+\frac{c_1}{u})^2+\frac{(1+\frac{c_2}{au}-\gamma_2\frac{c_2}{au}\frac{\ln^2u\Lambda}{u^2})^2}{1-2\gamma_2\frac{\ln^2u\Lambda}{u^2}+\gamma_2^2\frac{\ln^2u\Lambda}{u^2}}a_+^2}$$ for $j=3$.
	
	Consider first the case $j=1$. Then according to the Lemma $4.5$ from \cite{Vector} it is enough to prove that
	$$u^3e^{-\frac{\left(u+c_1\left(1-\frac{\ln^2u\Lambda}{u^2}\right)\right)^2+\left(u+\frac{c_2}{a}\left(1-\frac{\ln^2u\Lambda}{u^2}\right)\right)^2a_+^2}{2\left(1-\frac{\ln^2u\Lambda}{u^2}\right)}}=o(\pk{\vk B(1)-\vk c>\vk au}),$$
	which follows from the following equalities
	$$e^{-\frac{\left(u+c_1\left(1-\frac{\ln^2u\Lambda}{u^2}\right)\right)^2}{2\left(1-\frac{\ln^2u\Lambda}{u^2}\right)}}=o\left(u^{-k}e^{-\frac{(u+c_1)^2}{2}}\right)$$
	for $k>0$
	and
	$$e^{-\frac{\left(u+\frac{c_2}{a}\left(1-\frac{\ln^2u\Lambda}{u^2}\right)\right)^2a_+^2}{2\left(1-\frac{\ln^2u\Lambda}{u^2}\right)}}=O\left(e^{-\frac{(a_+u+c_2)^2}{2}}\right),$$
	as $u\to\IF$. Their proofs are too technical, hence we omit them here. 
	
	Let now $j=2$. As in the case $j=1$ it is sufficient to check that
	$$u^3e^{-\Bigg(\frac{(1+\frac{c_1}{u}-\gamma_1\frac{c_1}{u}\frac{\ln^2u\Lambda}{u^2})^2}{2(1-2\gamma_1\frac{\ln^2u\Lambda}{u^2}+\gamma_1^2\frac{\ln^2u\Lambda}{u^2})}+\frac{(1+\frac{c_2}{au})^2a_+^2}{2}\Bigg)u^2}=o(\pk{\vk B(1)>\vk au+\vk c}).$$
	One can see that
	$$e^{-\left(\frac{\left(1+\frac{c_1}{u}-\gamma_1\frac{c_1}{u}\frac{\ln^2u\Lambda}{u^2}\right)^2}{2\left(1-2\gamma_1\frac{\ln^2u\Lambda}{u^2}+\gamma_1^2\frac{\ln^2u\Lambda}{u^2}\right)}\right)u^2}=o\left(u^{-k}e^{-\frac{(u+c_1)^2}{2}}\right),$$
	as $u\to\IF$ for $k>0$ and
	$$e^{-\frac{(1+\frac{c_2}{au})^2a_+^2}{2}u^2}=O\left(e^{-\frac{(a_+u+c_2)^2}{2}}\right),$$
	as $u\to\IF$. The case $j=3$ can be checked in the similar way as the case $j=2$, so we omit it.
	
	Thus in all cases by the Lemma $4.5$ from \cite{Vector} we obtain that
	$$\pk{\exists_{(t,\vk s)\in A_1}:\vk B(t)-\vk ct-\vk\gamma(\vk B(\vk s)-\vk c\vk s)>\vk au}=o(\pk{\vk B(1)-\vk c>\vk au}).$$
	Consider now the case $a=0$. Observe that 
	$$\pk{\exists_{(t,\vk s)\in A_1}:\vk B(t)-\vk ct-\vk\gamma(\vk B(\vk s)-\vk c\vk s)>\vk au}\leq\pk{\exists_{(t,s)\in A_{1}^P}:B_1(t)- c_1t-\gamma_1(B_1(s_1)-c_1s_1)>u},$$
	where $A_{1}^P=P(A_{1})$ with a projection $P:\R^3\to\R^2$ such that $P((x_1,x_2,x_3))=(x_1,x_2)$.
	
	Moreover, for $a=0$
	$$\frac{1}{u}e^{-\frac{(u+c_1)^2}{2}}=O(\pk{\vk B(1)>\vk au+\vk c}).$$
	Then the proof is the same as for $a\neq 0$ but we consider the process $$B_1(t)-c_1t-\gamma_1(B_1(s_1)-c_1s_1)$$
	instead. Denote 
	$$A_{1,j}^P=P(A_{1,j}), j=1,2,3.$$
	Notice that 
	$$\pk{\exists_{(t,s_1)\in A_{1}^P}:B_1(t)- c_1t-\gamma_1(B_1(s)-c_1s_1)>u}=\pk{\exists_{(t,s_1)\in A_{1}^P}:X_u(t,s_1)>u},$$
	with 
	$$X_u(t,s_1)=\frac{B_1(1)-\gamma_1B_1(s_1)}{1+\frac{c_1t-\gamma_1c_1s_1}{u}}.$$
	Then following the same steps as in the case $a\neq 0$ one can check that
	$$\pk{\exists_{(t,s_1)\in A_{1}^P}:X_u(t,s_1)>u}=o\left(\frac{1}{u}e^{-\frac{(u+c_1)^2}{2}}\right)=o(\pk{\vk B(1)>\vk au+\vk c})$$
	for all $j=1,2,3$.
	Thus for the case $\gamma_1,\gamma_2>0$ the inequality (\ref{A_1}) holds.
	
	Assume now $\gamma_1=0$, $\gamma_2>0$.
	Then observe that since \(\vk X(t, \vk s)=\vk B(t)-\vk ct-\vk\gamma(\vk B(\vk s)-\vk c\vk s)\) is independent of $s_1$, we can see that
	$$\pk{\exists_{(t,\vk s)\in A_1}:\vk B(t)-\vk ct-\vk\gamma(\vk B(\vk s)-\vk c\vk s)>\vk au}=\pk{\exists_{(t,\vk s)\in A_1^0}:\vk B(t)-\vk ct-\vk\gamma^{0}(\vk B(\vk s)-\vk c\vk s)>\vk au},$$
	where 
	$$A_1^0=A_1\cap\{(x_1,x_2,x_3)\in\R^3:x_2=0\}\text{ and }\vk\gamma^{0}=\left(\frac{1}{2},\gamma_2\right).$$
	
	Since $\vk\gamma^{0}$ has only positive components now and $A_1^0\subset A_1$, by the discussed before case $\gamma_1,\gamma_2>0$ we have that
	\begin{align*}
		\pk{\exists_{(t,\vk s)\in A_1^0}:\vk B(t)-\vk ct-\vk\gamma^{0}(\vk B(\vk s)-\vk c\vk s)>\vk au}&\leq\pk{\exists_{(t,\vk s)\in A_1}:\vk B(t)-\vk ct-\vk\gamma^{0}(\vk B(\vk s)-\vk c\vk s)>\vk au}\\&=o(\pk{\vk B(1)>\vk au+\vk c}).
	\end{align*}
	The cases $\gamma_1>0,\gamma_2=0$ and $\gamma_1=0,\gamma_2=0$ can be done similarly. 
	
	If $\gamma_1>0,\gamma_2=0$, then take 
	$$A_1^0=A_1\cap\{(x_1,x_2,x_3)\in\R^3:x_3=0\} \text{ and } \vk\gamma^{0}=\left(\gamma_1, \frac{1}{2}\right)$$ 
	
	and if $\gamma_1=0,\gamma_2=0$, then 
	$$A_1^0=A_1\cap\{(x_1,x_2,x_3)\in\R^3:x_2=x_3=0\} \text{ and } \vk\gamma^{0}=\left(\frac{1}{2}, \frac{1}{2}\right).$$
	
	In both cases we obtain the same chain of inequalities as for the case $\gamma_1=0$, $\gamma_2>0$.

	\subsection{Proof of (\ref{A_0})}
	Assume first $\gamma_1,\gamma_2>0$. Let now $a>0$.
	
	Denote  
	$$\Delta_{k,l_1,l_2}:=\left[1-\frac{(k+1)\Lambda}{u^2}, 1-\frac{k\Lambda}{u^2}\right]\times\left[\frac{l_1\Lambda}{u^2}, \frac{(l_1+1)\Lambda}{u^2}\right]\times\left[\frac{l_2\Lambda}{u^2}, \frac{(l_2+1)\Lambda}{u^2}\right],$$
	where $k,l_1,l_2$ are non-negative integers such that $k^2+l_1^2+l_2^2>0$ and $k,l_1,l_2$ are not greater than $\ln^2u$.
	
	The idea of the proof here is similar to the proof of Lemma \ref{M(u,Lambda)}.
	Fix numbers $k,l_1,l_2$ and let 
	$$\vk X(t,\vk s)=\vk B(t)-\vk\gamma\vk B(\vk s).$$
	Then our goal it to estimate the following probability
	$$\pk{\exists_{(t,\vk s)\in \Delta_{k,l_1,l_2}}:\vk X(t,\vk s)-\vk ct+\vk\gamma\vk c\vk s>\vk au}.$$
	Denote by $\phi_{k,l_1,l_2}$ the density function of the vector $\vk X\left(1-\frac{k\Lambda}{u^2},\frac{l_1\Lambda}{u^2},\frac{l_2\Lambda}{u^2}\right).$
	
	Notice that the covariance function of $\vk X\left(1-\frac{k\Lambda}{u^2},\frac{l_1\Lambda}{u^2},\frac{l_2\Lambda}{u^2}\right)$ is equal to
	$$\Sigma_{k,l_1,l_2}=\begin{pmatrix}1-\frac{k\Lambda}{u^2}+\gamma_1^2\frac{l_1\Lambda}{u^2}-2\gamma_1\frac{l_1\Lambda}{u^2} & 0\\0 & 1-\frac{k\Lambda}{u^2}+\gamma_2^2\frac{l_2\Lambda}{u^2}-2\gamma_2\frac{l_2\Lambda}{u^2}\end{pmatrix}$$
	and, hence,
	$$\Sigma_{k,l_1,l_2}^{-1}=\begin{pmatrix}\frac{1}{1-\frac{k\Lambda}{u^2}+\gamma_1^2\frac{l_1\Lambda}{u^2}-2\gamma_1\frac{l_1\Lambda}{u^2}} & 0\\0 & \frac{1}{1-\frac{k\Lambda}{u^2}+\gamma_2^2\frac{l_2\Lambda}{u^2}-2\gamma_2\frac{l_2\Lambda}{u^2}}\end{pmatrix}.$$
	Thus 
	$$\phi_{k,l_1,l_2}(x_1,x_2)=\frac{1}{2\pi\sqrt{\left(1-\frac{(k+2\gamma_1l_1-\gamma_1^2l_1)\Lambda}{u^2}\right)\left(1-\frac{(k+2\gamma_2l_2-\gamma_2^2l_2)\Lambda}{u^2}\right)}}e^{-\left(\frac{x_1^2}{2\left(1-\frac{(k+2\gamma_1l_1-\gamma_1^2l_1)\Lambda}{u^2}\right)}+\frac{x_2^2}
		{2\left(1-\frac{(k+2\gamma_2l_2-\gamma_2^2l_2)\Lambda}{u^2}\right)}\right)}.$$
	Then one can check that for $u$ large enough 
	$$\phi_{k,l_1,l_2}\left(\vk au+\vk c-\frac{\vk x}{u}\right)\leq 2\phi(\vk au+\vk c)e^{x_1\left(1+\frac{sign(x_1)}{2}\right)+ax_2\left(1+\frac{sign(x_2)}{2}\right)}e^{-\frac{(1+a^2)k\Lambda+(2\gamma_1-\gamma_1^2)l_1\Lambda+(2\gamma_2a^2-\gamma_2^2a^2)l_2\Lambda}{2}},$$
	where $\phi$ is the density function of the vector $\vk B(1)$.

	To shorten further designations denote
	$$t_{uk}=1-\frac{k\Lambda+t}{u^2},\vk s_{ul}=(s_{ul,1},s_{ul,2})=\left(\frac{l_1\Lambda+s_1}{u^2},\frac{l_2\Lambda+s_2}{u^2}\right),$$ $$r_{uk}=1-\frac{k\Lambda}{u^2}, \vk q_{ul}=(q_{ul,1},q_{ul,2})=\left(\frac{l_1\Lambda}{u^2},\frac{l_2\Lambda}{u^2}\right)$$ 
	for $(t,s_1,s_2)\in [0,\Lambda]^3$ and as before 
	$$\vk u_{x}=\vk au+\vk c-\frac{\vk x}{u}.$$
    Conditioning, we obtain
	\begin{align*}
		\frac{1}{u^2}\int_{\R^2}\pk{\exists_{(t,\vk s)\in [0,\Lambda]^3}:\vk X(t_{uk},\vk s_{ul})-\vk ct_{uk}+\vk\gamma\vk c\vk s_{ul}>\vk au\bigg|\vk X(r_{uk},\vk q_{ul})=\vk u_{x}}\phi_{k,l_1,l_2}(\vk u_{x})d\vk x.
	\end{align*}
	
	Observe that $\vk X(t_{uk},\vk s_{ul})$ can be rewritten in the following view
	$$\vk X(t_{uk},\vk s_{ul})=cov(\vk X(t_{uk},\vk s_{ul}),\vk X(r_{uk},\vk q_{ul}))\Sigma_{k,l_1,l_2}^{-1}\vk X(r_{uk},\vk q_{ul})+\vk D_{u,k,l}(t,\vk s),$$
	where $\vk D_{u,k,l}(t,\vk s)$ and $\vk X(r_{uk},\vk q_{ul})$ are independent processes and $\Sigma_{k,l_1,l_2}$ is a covariance matrix for the vector $\vk X(r_{uk},\vk q_{ul})$.
	Thus the event 
	$$\left\{\exists_{(t,\vk s)\in [0,\Lambda]^3}:\vk X(t_{uk},\vk s_{ul})-\vk ct_{uk}+\vk\gamma\vk c\vk s_{ul}>\vk au\bigg|\vk X(r_{uk},\vk q_{ul})=\vk u_{x}\right\}$$
	can be rewritten as follows
	$$\left\{\exists_{(t,\vk s)\in [0,\Lambda]^3}:cov(\vk X(t_{uk},\vk s_{ul}),\vk X(r_{uk},\vk q_{ul}))\Sigma_{k,l_1,l_2}^{-1}\vk u_{x}+\vk D_{u,k,l}(t,\vk s)-\vk ct_{uk}+\vk\gamma\vk c\vk s_{ul}>\vk au\right\}.$$
	Consider now $cov(\vk X(t_{uk},\vk s_{ul}),\vk X(r_{uk},\vk q_{ul}))$.
	Since the Brownian motions $B_1$ and $B_2$ are independent, the covariance matrix $cov(\vk X(t_{uk},\vk s_{ul}),\vk X(r_{uk},\vk q_{ul}))$ is equal to
	$$\begin{pmatrix}cov(X_1(t_{uk},\vk s_{ul}),X_1(r_{uk},\vk q_{ul})) & 0\\0 & cov(X_2(t_{uk},\vk s_{ul}), X_2(r_{uk},\vk q_{ul}))\end{pmatrix}.$$
	Calculate $cov(X_i(t_{uk},\vk s_{ul}),X_i(r_{uk},\vk q_{ul}))$. It is equal to
	$$cov(B_i(t_{uk})-\gamma_iB_i(s_{ul,i}),B_i(r_{uk})-\gamma_iB_i(q_{ul,i}))=t_{uk}-\gamma_is_{ul,i}-\gamma_iq_{ul,i}+\gamma_i^2q_{ul,i}.$$
	Thus $$cov(\vk X(t_{uk},\vk s_{ul}),\vk X(r_{uk},\vk q_{ul}))=\begin{pmatrix}t_{uk}-\gamma_1s_{ul,1}-\gamma_1q_{ul,1}+\gamma_1^2q_{ul,1} & 0\\0 & t_{uk}-\gamma_2s_{ul,2}-\gamma_2q_{ul,2}+\gamma_2^2q_{ul,2}\end{pmatrix}.$$
	Since
	$$\Sigma_{k,l_1,l_2}=\begin{pmatrix}r_{uk}-2\gamma_1q_{ul,1}+\gamma_1^2q_{ul,1} & 0\\0 & r_{uk}-2\gamma_2q_{ul,2}+\gamma_2^2q_{ul,2}\end{pmatrix}$$
	we conclude that
	\begin{align*}
		cov(\vk X(t_{uk},\vk s_{ul}),\vk X(r_{uk},\vk q_{ul}))\Sigma_{k,l_1,l_2}^{-1}&=\begin{pmatrix}\frac{t_{uk}-\gamma_1s_{ul,1}-\gamma_1q_{ul,1}+\gamma_1^2q_{ul,1}}{r_{uk}-2\gamma_1q_{ul,1}+\gamma_1^2q_{ul,1}} & 0\\0 & \frac{t_{uk}-\gamma_2q_{ul,2}-\gamma_2q_{ul,2}+\gamma_2^2q_{ul,2}}{r_{uk}-2\gamma_2q_{ul,2}+\gamma_2^2q_{ul,2}}\end{pmatrix}\\&=\begin{pmatrix}1-\frac{t+\gamma_1s_1}{u^2-k\Lambda-2\gamma_1l_1\Lambda+\gamma_1^2l_1\Lambda} & 0\\0 & 1-\frac{t+\gamma_2s_2}{u^2-k\Lambda-2\gamma_2l_2\Lambda+\gamma_2^2l_2\Lambda}\end{pmatrix}.
	\end{align*}
	Now the idea of the proof is the same as in the proof of Lemma \ref{M(u,Lambda)}, when we justified the use of dominated convergence theorem. 
	
	Let first $x_1,x_2>0$. Then for large enough $u$
	\begin{align*}
		&\pk{\exists_{(t,\vk s)\in [0,\Lambda]^3}:cov\left(\vk X\left(t_{uk},\vk s_{ul}\right),\vk X\left(r_{uk},\vk q_{ul}\right)\right)\Sigma_{k,l_1,l_2}^{-1}\vk u_{x}+\vk D_{u,k,l}(t,\vk s)-\vk ct_{uk}+\vk\gamma\vk c\vk s_{ul}>\vk au}&
		\\&\leq \pk{\exists_{(t,\vk s)\in [0,\Lambda]^3}:\sum_{i=1}^{2}u\left(cov\left(\vk X(t_{uk},\vk s_{ul}),\vk X(r_{uk},\vk q_{ul})\right)\Sigma_{k,l_1,l_2}^{-1}\vk u_{x}+\vk D_{u,k,l}(t,\vk s)-\vk ct_{uk}+\vk\gamma\vk c\vk s_{ul}\right)_i>(1+a)u^2}&
		\\&\leq \pk{\exists_{(t,\vk s)\in [0,\Lambda]^3}:\sum_{i=1}^{2}u\left(\left(1-\frac{t+\gamma_is_i}{u^2-k\Lambda-2\gamma_il_i\Lambda+\gamma_i^2l_i\Lambda}\right)\vk u_{x}+\vk D_{u,k,l}(t,\vk s)-\vk ct_{uk}+\vk\gamma\vk c\vk s_{ul}\right)_i>(1+a)u^2}&
		\\&\leq\pk{\exists_{(t,\vk s)\in [0,\Lambda]^3}:u\sum_{i=1}^{2}(D_{u,k,l}(t,\vk s))_i>\frac{x_1+x_2}{2}}.&
	\end{align*}
	Here we used that $k,l_1,l_2=o(u)$ and $t,s_1,s_2$ are bounded.
	
	Now our aim is to use Borel's inequality for the one-dimensional process
	$u\sum_{i=1}^{2}(D_{u,k,l}(t,\vk s))_i$ and we need to check that its variance is uniformly bounded by some constant, which does not depend on $u,k,l_1,l_2$.
	First notice since $B_1$ and $B_2$ are independent processes,
	$$Var\left(u\sum_{i=1}^{2}(D_{u,k,l}(t,\vk s))_i\right)=u^2\sum_{i=1}^{2}Var((D_{u,k,l}(t,\vk s))_i).$$
	Thus it is sufficient to estimate uniformly $u^2Var((D_{u,k,l}(t,\vk s))_i)$ for $i=1$ and $i=2$.
	Without loss of generality assume that $i=1$. Then
	\begin{align*}
		&u^2Var((D_{u,k,l}(t,\vk s))_1)=u^2Var\Bigg(B_1(t_{uk})-\gamma_1B_1(s_{ul,1})&
		\\&-\left(1-\frac{t+\gamma_1s_1}{u^2-k\Lambda-2\gamma_1l_1\Lambda+\gamma_1^2l_1\Lambda}\right)\left(B_1\left(r_{uk}\right)-\gamma_1B_1\left(q_{ul,1}\right)\right)\Bigg)&
		\\&\leq 2u^2Var\left(B_1(t_{uk})-\left(1-\frac{t+\gamma_1s_1}{u^2-k\Lambda-2\gamma_1l_1\Lambda+\gamma_1^2l_1\Lambda}\right)B_1\left(r_{uk}\right)\right)&
		\\&+2u^2Var\left(\gamma_1B_1(s_{ul,1})-\left(1-\frac{t+\gamma_1s_1}{u^2-k\Lambda-2\gamma_1l_1\Lambda+\gamma_1^2l_1\Lambda}\right)\gamma_1B_1\left(q_{ul,1}\right)\right)&.
	\end{align*}
	Estimate now $$2u^2Var\left(\gamma_1B_1(s_{ul,1})-\left(1-\frac{t+\gamma_1s_1}{u^2-k\Lambda-2\gamma_1l_1\Lambda+\gamma_1^2l_1\Lambda}\right)\gamma_1B_1\left(q_{ul,1}\right)\right).$$ 
	The other term from the last sum can be estimated in the similar way. Now,
	\begin{align*}
		&2u^2Var\left(\gamma_1B_1(s_{ul,1})-\left(1-\frac{t+\gamma_1s_1}{u^2-k\Lambda-2\gamma_1l_1\Lambda+\gamma_1^2l_1\Lambda}\right)\gamma_1B_1\left(q_{ul,1}\right)\right)&
		\\&=2u^2\gamma_1^2Var\left(B_1(s_{ul,1})-\left(1-\frac{t+\gamma_1s_1}{u^2-k\Lambda-2\gamma_1l_1\Lambda+\gamma_1^2l_1\Lambda}\right)B_1\left(q_{ul,1}\right)\right)&
		\\&=2u^2\gamma_1^2\left(s_{ul,1}-2\left(1-\frac{t+\gamma_1s_1}{u^2-k\Lambda-2\gamma_1l_1\Lambda+\gamma_1^2l_1\Lambda}\right)q_{ul,1}+q_{ul,1}\left(1-\frac{t+\gamma_1s_1}{u^2-k\Lambda-2\gamma_1l_1\Lambda+\gamma_1^2l_1\Lambda}\right)^2\right)&
		\\&\leq2u^2\gamma_1^2\left(s_{ul,1}-2\left(1-\frac{t+\gamma_1s_1}{u^2-k\Lambda-2\gamma_1l_1\Lambda+\gamma_1^2l_1\Lambda}\right)q_{ul,1}+q_{ul,1}\left(1-\frac{t+\gamma_1s_1}{u^2-k\Lambda-2\gamma_1l_1\Lambda+\gamma_1^2l_1\Lambda}\right)\right)&
		\\&\leq2u^2\gamma_1^2\left(s_{ul,1}-\left(1-\frac{t+\gamma_1s_1}{u^2-k\Lambda-2\gamma_1l_1\Lambda+\gamma_1^2l_1\Lambda}\right)q_{ul,1}\right)&
		\\&=2\gamma_1^2s_1+2\gamma_1^2\frac{t+\gamma_1s_1}{u^2-k\Lambda-2\gamma_1l_1\Lambda+\gamma_1^2l_1\Lambda}l_1\Lambda<4\gamma_1^2\Lambda.&
	\end{align*}
	Thus there is a constant $C'_1$, which does not depend on the choice of $k,l_1$, such that
	$$u^2Var((D_{u,k,l}(t,\vk s))_1)\leq C'_1.$$
	In a similar way one can check that there exists a constant $C'_2$, which does not depend on the choice of $k,l_2$, such that
	$$u^2Var((D_{u,k,l}(t,\vk s))_2)\leq C'_2.$$
	
	Thus using \cite[Lemma 4.5]{Vector} we obtain that
	\begin{align*}
		&\pk{\exists_{(t,\vk s)\in [0,\Lambda]^3}:cov(\vk X(t_{uk},\vk s_{ul}),\vk X(r_{uk},\vk q_{ul}))\Sigma_{k,l_1,l_2}^{-1}\vk u_{x}+\vk D_{u,k,l}(t,\vk s)-\vk ct_{uk}+\vk\gamma\vk c\vk s_{ul}>\vk au}&
		\\&\leq C_0(e^{-C_1(x_1+x_2)^2}\mathbb{I}(x_1+x_2\geq C_2)+\mathbb{I}(x_1+x_2<C_2))\leq C_0(e^{-C_1(x_1^2+x_2^2)}\mathbb{I}(x_1+x_2\geq C_2)+\mathbb{I}(x_1+x_2<C_2))&
	\end{align*}
	for some positive constants $C_0,C_1$ and $C_2$.
	
	Let now $x_1\leq 0, x_2>0$. Then
	\begin{align*}
		&\pk{\exists_{(t,\vk s)\in [0,\Lambda]^3}:cov(\vk X(t_{uk},\vk s_{ul}),\vk X(r_{uk},\vk q_{ul}))\Sigma_{k,l_1,l_2}^{-1}\vk u_{x}+\vk D_{u,k,l}(t,\vk s)-\vk ct_{uk}+\vk\gamma\vk c\vk s_{ul}>\vk au}&
		\\&\leq \pk{\exists_{(t,\vk s)\in [0,\Lambda]^3}:u(cov(\vk X(t_{uk},\vk s_{ul}),\vk X(r_{uk},\vk q_{ul}))\Sigma_{k,l_1,l_2}^{-1}\vk u_{x}+\vk D_{u,k,l}(t,\vk s)-\vk ct_{uk}+\vk\gamma\vk c\vk s_{ul})_1>u^2}&
		\\&\leq E_0(e^{-E_1x_1^2}\mathbb{I}(x_1\geq E_2)+\mathbb{I}(x_1<E_2))\leq C_0(e^{-E_1x_1^2}\mathbb{I}(x_1\geq E_2)+\mathbb{I}(x_1<E_2))&
	\end{align*}
	for some constants $E_0,E_1,E_2>0$. 
	
	If $x_1>0, x_2\leq 0$. Then by the similar reasons we obtain that 
	\begin{align*}
		&\pk{\exists_{(t,\vk s)\in [0,\Lambda]^3}:cov(\vk X(t_{uk},\vk s_{ul}),\vk X(r_{uk},\vk q_{ul}))\Sigma_{k,l_1,l_2}^{-1}\vk u_{x}+\vk D_{u,k,l}(t,\vk s)-\vk ct_{uk}+\vk\gamma\vk c\vk s_{ul}>\vk au}&
		\\&\leq R_0(e^{-R_1x_2^2}\mathbb{I}(x_2\geq R_2)+\mathbb{I}(x_2<R_2))\leq R_0(e^{-R_1x_2^2}\mathbb{I}(x_2\geq R_2)+\mathbb{I}(x_2<R_2))&
	\end{align*}
	for some positive constants $R_0,R_1$ and $R_2$.
	
	To sum up, according to \cite[Lemma 2]{Lemma2}, as $u\to\IF$, we obtain the following chain of inequalities
	\begin{align*}
		&\pk{\exists_{(t,\vk s)\in \Delta_{k,l_1,l_2}}:\vk X(t,\vk s)-\vk ct+\vk\gamma\vk c\vk s>\vk au}\leq\frac{2}{u^2}\phi(\vk au+\vk c)e^{-\frac{(1+a^2)k\Lambda+(2\gamma_1-\gamma_1^2)l_1\Lambda+(2\gamma_2a-\gamma_2^2a^2)l_2\Lambda}{2}}\int_{\R^2 }h(\vk x)d\vk x&
		\\&\leq 4\pk{\vk B(1)>\vk au+\vk c}e^{-\frac{(1+a^2)k\Lambda+(2\gamma_1-\gamma_1^2)l_1\Lambda+(2\gamma_2a^2-\gamma_2^2a^2)l_2\Lambda}{2}}\int_{\R^2 }h(\vk x)d\vk x,&
	\end{align*}
	where 
	$$h(x_1,x_2):=\begin{cases}e^{\frac{x_1+ax_2}{2}}, \ \ \ x_1,x_2\leq 0;\\E_0(e^{-E_1x_2^2}\mathbb{I}(x_2\geq E_2)+\mathbb{I}(x_2<E_2))e^{2ax_2}, \ \ \ x_1\leq 0,x_2>0;\\R_0(e^{-R_1x_1^2}\mathbb{I}(x_1\geq R_2)+\mathbb{I}(x_1<R_2))e^{2x_1}, \ \ \ x_1> 0,x_2\leq 0;\\C_0(e^{-C_1(x_1^2+x_2^2)}\mathbb{I}(x_1+x_2\geq C_2)+\mathbb{I}(x_1+x_2<C_2))e^{2x_1+2ax_2}, \ \ \ x_1> 0,x_2>0.\end{cases}$$
	It follows that $h(x_1,x_2)$ is integrable on $\R^2$. Moreover, observe that we chose the function independent of parameters $k,l_1,l_2$. Thus since \(\gamma_1,\gamma_2\in (0,2)\), for all non-negative integers $k,l_1,l_2$ such that $k^2+l_1^2+l_2^2\neq 0$ there exist constants $C^{*},C_0^{*},C_1^{*},C_2^{*}$, which are not dependent on $k,l_1,l_2$, such that
	$$\pk{\exists_{(t,\vk s)\in\Delta_{k,l_1,l_2}}:\vk B(t)-\vk ct-\vk\gamma(\vk B(\vk s)-\vk c\vk s)>\vk au}\leq C^{*}e^{-(C_0^{*}k+C_1^{*}l_1+C_2^{*}l_2)\Lambda}\pk{\vk B(1)>\vk au+\vk c}.$$
	For $\Lambda>1$ (we are interested only in large $\Lambda$)
	\begin{align*}
		&\frac{\pk{\exists_{(t,\vk s)\in A_0}:\vk B(t)-\vk c(t)-\vk\gamma(\vk B(\vk s)-\vk c\vk s)>\vk au}}{\pk{\vk B(1)>\vk au+\vk c}}\leq\sum_{k=0}^{[\ln^2u]}\sum_{l_1=0}^{[\ln^2u]}\sum_{l_2=0}^{[\ln^2u]}C^{*}e^{-(C_0^{*}k+C_1^{*}l_1+C_2^{*}l_2)\Lambda}-C^{*}&
		\\&\leq C^{*}\Bigl(\sum_{k=0}^{\IF}e^{-C_0^{*}k\Lambda}\Bigr)\Bigl(\sum_{l_1=0}^{\IF}e^{-C_1^{*}l_1\Lambda}\Bigr)\Bigl(\sum_{l_2=0}^{\IF}e^{-C_2^{*}l_2\Lambda}\Bigr)-C^{*}&
		\\&\leq C^{*}\Bigl(1+e^{-C_0^{*}\Lambda}\frac{1}{1-e^{-C_0^{*}\Lambda}}\Bigr)\Bigl(1+e^{-C_1^{*}\Lambda}\frac{1}{1-e^{-C_1^{*}\Lambda}}\Bigr)\Bigl(1+e^{-C_2^{*}\Lambda}\frac{1}{1-e^{-C_2^{*}\Lambda}}\Bigr)-C^{*}&
		\\&\leq C^{*}\Bigl(1+e^{-C_0^{*}\Lambda}\frac{1}{1-e^{-C_0^{*}}}\Bigr)\Bigl(1+e^{-C_1^{*}\Lambda}\frac{1}{1-e^{-C_1^{*}}}\Bigr)\Bigl(1+e^{-C_2^{*}\Lambda}\frac{1}{1-e^{-C_2^{*}}}\Bigr)-C^{*}&
		\\&\leq G_0(e^{-C_0^{*}\Lambda}+e^{-C_1^{*}\Lambda}+e^{-C_2^{*}\Lambda})\leq 3G_0e^{-\min(C_0^{*},C_1^{*},C_2^{*})\Lambda}=:G_1e^{-G_2\Lambda},&
	\end{align*}
	as $u\to\IF$, where $G_0, G_1, G_2$ are some positive constants, which does not dependent on $\Lambda$.
	
	Consider now the case $a\leq 0$. The idea of the proof here is exactly the same as in the case $a>0$, so we discuss it briefly. Recall that for $a\leq 0$
	$$\pk{\vk B(1)>\vk au+\vk c}\sim F_0\frac{1}{u}e^{\frac{-(u+c_1)^2}{2}}$$
	for some constant $F_0$. Thus there exists a constant $F_1$ such that for large $u$ 
	$$\pk{B_1(1)>u+c_1}\leq F_1\pk{\vk B(1)>\vk au+\vk c}.$$
	On the other hand
	\begin{align*}
		&\pk{\exists_{(t,\vk s)\in A_0}:\vk B(t)-\vk ct-\vk\gamma(\vk B(\vk s)-\vk c\vk s)>\vk au}&
		\\&\leq\pk{\exists_{(t,s_1)\in \left(\left[1-\frac{\ln^2u\Lambda}{u^2}\right]\times\left[0,\frac{\ln^2u\Lambda}{u^2}\right]\right)\setminus \left(\left[1-\frac{\Lambda}{u^2}\right]\times\left[0,\frac{\Lambda}{u^2}\right]\right)}:B_1(t)-c_1t-\gamma_1(B_1(s_1)-c_1s_1)>u}.&
	\end{align*}
	Thus it is enough to find positive constants $F, F^{*}>0$ such that
	\begin{align*}
		\pk{\exists_{(t,s_1)\in \left(\left[1-\frac{\ln^2u\Lambda}{u^2}\right]\times\left[0,\frac{\ln^2u\Lambda}{u^2}\right]\right)\setminus \left(\left[1-\frac{\Lambda}{u^2}\right]\times\left[0,\frac{\Lambda}{u^2}\right]\right)}:B_1(t)-c_1t-\gamma_1(B_1(s_1)-c_1s_1)>u}\\\leq F^{*}e^{-F\Lambda}\pk{B_1(1)>u+c_1}.
	\end{align*}
	We divide now the set
	$$\left(\left[1-\frac{\ln^2u\Lambda}{u^2},1\right]\times\left[0,\frac{\ln^2u\Lambda}{u^2}\right]\right)\setminus \left(\left[1-\frac{\Lambda}{u^2},1\right]\times\left[0,\frac{\Lambda}{u^2}\right]\right)$$ into rectangles $$\Delta_{k,l}:=\left[1-\frac{(k+1)\Lambda}{u^2},1-\frac{k\Lambda}{u^2}\right]\times\left[\frac{l\Lambda}{u^2},\frac{(l+1)\Lambda}{u^2}\right]$$ for $k^2+l^2\neq 0$. 
	Using similar ideas as in the case $a>0$ we find positive constants $C^{*},C_0^{*},C_1^{*}$ such that
	$$\pk{\exists_{(t,s_1)\in\Delta_{k,l}}:B_1(t)-c_1t-\gamma_1(B_1(s_1)-c_1s_1)>u}\leq C^{*}e^{-(C_0^{*}k+C_1^{*}l_1)\Lambda}\pk{\vk B(1)>\vk au+\vk c}$$
	and then by summing estimates for all sets $\Delta_{k,l}$ with $k^2+l^2\neq 0$
	we can find constants $G_0, G_1,G_2>0$ such that
	$$\pk{\exists_{(t,\vk s)\in A_0}:\vk B(t)-\vk ct-\vk\gamma(\vk B(\vk s)-\vk c\vk s)>\vk au}\leq G_0e^{-G_2\Lambda}\pk{B_1(1)>u+c_1}\leq G_1e^{-G_2\Lambda}\pk{\vk B(1)>\vk au+\vk c}.$$
	That finishes the proof for the case $\gamma_1,\gamma_2>0$.
	
	The cases $\gamma_1=0$ or $\gamma_2=0$ can be derived from the case $\gamma_1,\gamma_2>0$ using the same argument as in the proof of inequality (\ref{A_1}).
	\section*{Acknowledgments}
	Partial  supported by SNSF Grant 200021-196888 is kindly acknowledged.
	
	\def\polhk#1{\setbox0=\hbox{#1}{\ooalign{\hidewidth
				\lower1.5ex\hbox{`}\hidewidth\crcr\unhbox0}}}

\end{document}